\documentclass[final]{siamltex}
\usepackage{amstext}
\usepackage{graphicx,subfigure,afterpage}
\usepackage{amsfonts}
\usepackage{hyperref}
\usepackage{amsmath}
\usepackage{color,fancybox,shadow}

\usepackage{url}

\makeatletter
\def\url@leostyle{%
  \@ifundefined{selectfont}{\def\UrlFont{\sf}}{\def\UrlFont{\small\ttfamily}}}
\makeatother
\newcommand{\beq}{\begin{equation}}
\newcommand{\eeq}{\end{equation}}

\long\def\drop#1{}

\def\R{\mathbb R}

\def\L{\mathcal L}

\def\Xint#1{\mathchoice
   {\XXint\displaystyle\textstyle{#1}}%
   {\XXint\textstyle\scriptstyle{#1}}%
   {\XXint\scriptstyle\scriptscriptstyle{#1}}%
   {\XXint\scriptscriptstyle\scriptscriptstyle{#1}}%
   \!\int}
\def\XXint#1#2#3{{\setbox0=\hbox{$#1{#2#3}{\int}$}
     \vcenter{\hbox{$#2#3$}}\kern-.5\wd0}}

\def\dashint{\Xint-}
\begin{document}

\author{Rustum Choksi\thanks{Dept.~of Mathematics,  McGill University, Montreal, Canada, {\tt rchoksi@math.mcgill.ca}} 
\and Mirjana Maras\thanks{Dept.~of Mathematics, Simon Fraser University, Burnaby, Canada, {\tt mmaras@sfu.ca}}
\and J.F. Williams\thanks{Dept.~of Mathematics, Simon Fraser University, Burnaby, Canada, {\tt jfw@math.sfu.ca}}}

\title{2D Phase Diagram for Minimizers of a  Cahn-Hilliard Functional with Long-range Interactions}

\maketitle

\begin{abstract} 
This paper presents a detailed asymptotic and numerical investigation of the phase diagram for global minimizers to a {\it Cahn-Hilliard functional with long-range interactions} in two space dimensions. We introduce a small parameter measuring perturbation from the minimal order-disorder transition, and 
derive asymptotic estimates for stability regions as the parameter  tends to zero.  
Based upon the $H^{-1}$ gradient flow, we introduce a hybrid numerical method  to navigate through the complex  energy landscape and  access the ground state of the functional. We use this method to numerically compute the phase diagram. 
Our asymptotic predictions  show  surprisingly good agreement with our numerical results.
\end{abstract} 


\section{Introduction} 
\begin{figure}[t!]
\centerline{{\includegraphics[width=.21\textwidth]{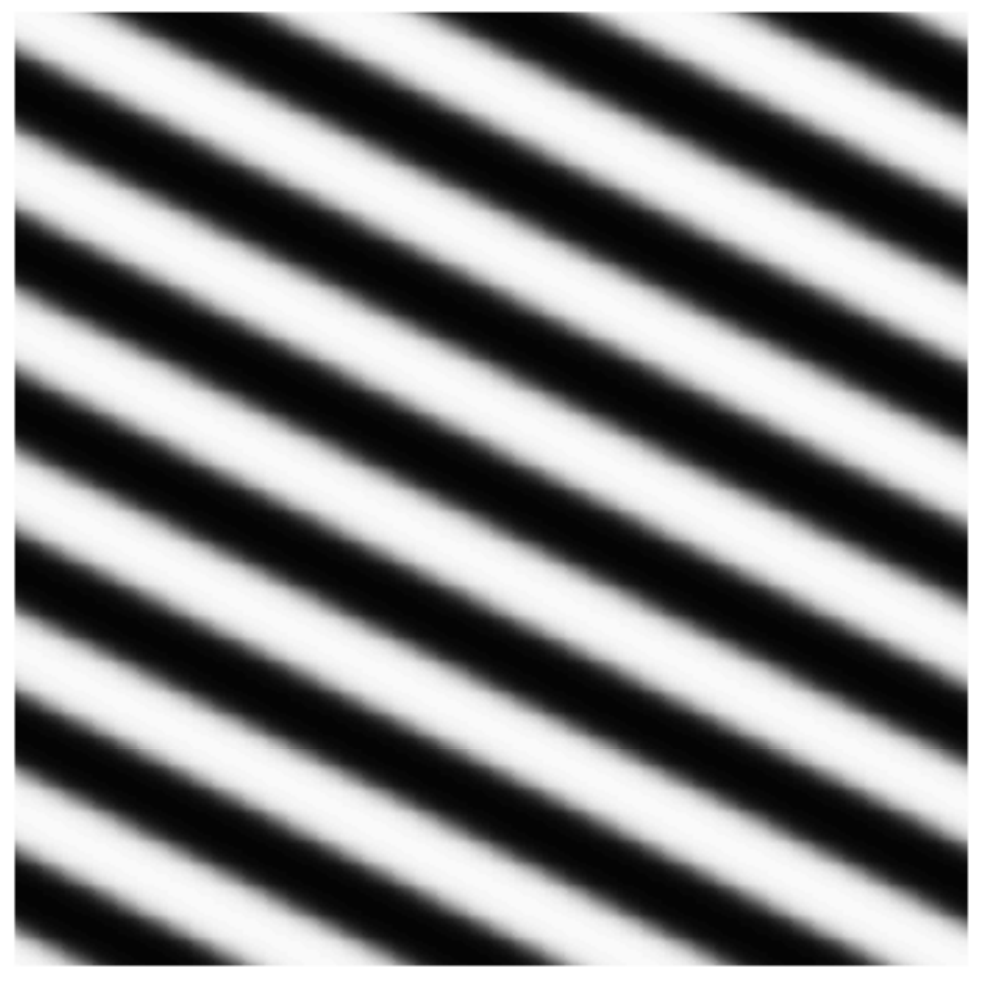}}
{\includegraphics[width=.21\textwidth]{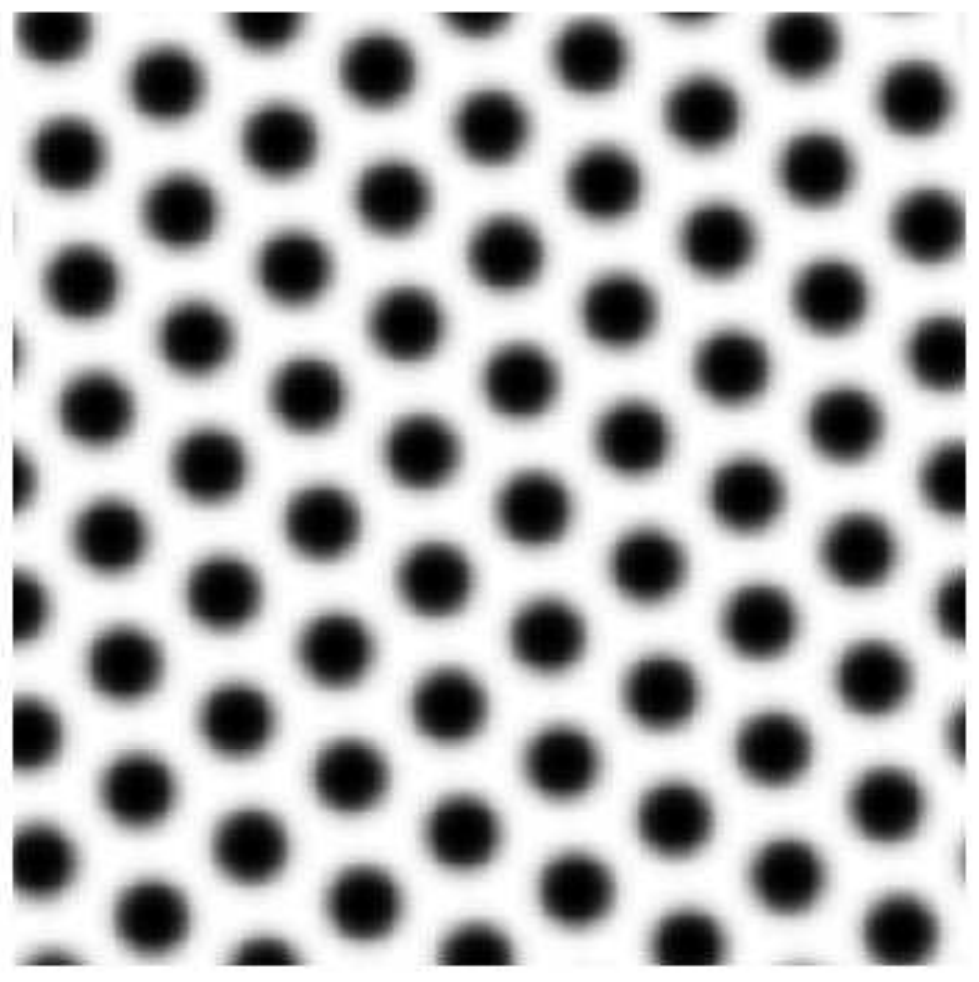}}
{\includegraphics[width=.21\textwidth]{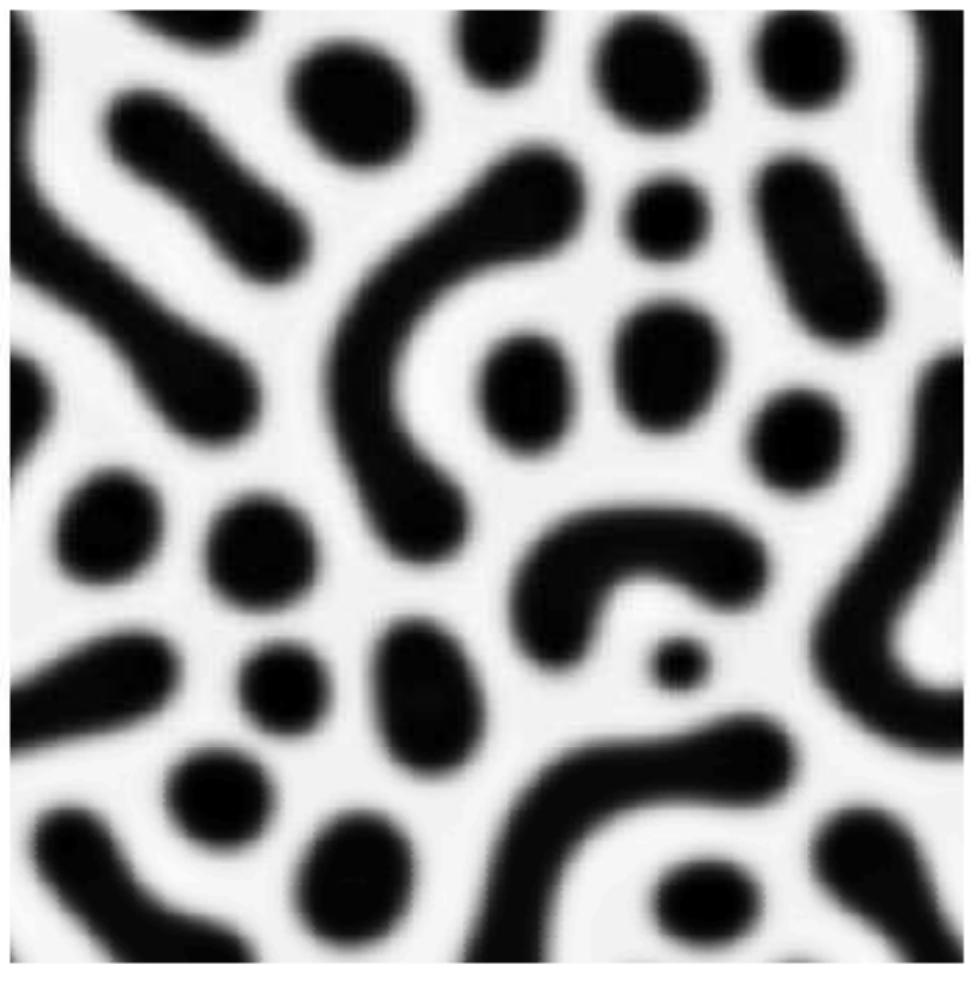}}}
\caption{Typical long-time solutions to (\ref{pde}) with $\gamma = 10$. Left: Lamellae, $m=0$. This is a global minimizing state.
Center: Hexagonally packed spots, $m=.3$. This is a global minimizing state. Right: A mixed state, $m=.4$. This is a typical metastable solution. All figures were computed on a domain of size $4\pi\times4\pi$. Here and in the following figures $u = 1$ is represented in black and $u = -1$ in white.}
\label{Fig:Intro}
\end{figure}

In this article, we asymptotically and numerically address the phase diagram with respect to the parameters $\gamma$ and $m$ for the following mass-constrained  
variational problem: For $\gamma >0$ and $m \in (-1,1)$, minimize 
\begin{equation}\label{OK}   \int_{\Omega} \Bigl(  \frac{1}{\gamma^2} |\nabla u|^2  +
\frac{(1 - u^2)^2}{4}\Bigr) \, d {x}
 \, + \, 
  \int_{\Omega}  \int_{\Omega} \, G({ x}, { y})  \left(u( x)  - {m}  \right)
\left(u( y)  - m \right) \,    d { x} \, d{y}, \end{equation}
over all $u$ with $\dashint_\Omega u \, dx= m$. 
Here  $G$ denotes the Green's function of $-\triangle$ on a  cubic domain $\Omega: = [0,L]^n \subset \R^n$ with periodic boundary conditions.  We refer to functional (\ref{OK}) as a Cahn-Hilliard functional ({\it cf.} \cite{CH}) with long-range interactions\footnote{This functional is also commonly referred 
as a  {\it Ginzburg-Landau} functional with {\it competing or Coulomb-type} interactions \cite{SD, Mu}, or because of its relation to self-assembly of diblock copolymers, the Ohta-Kawasaki functional \cite{OK} or simply as {\it the diblock copolymer problem} \cite{RW}. 
Our reference and labeling of {\it Cahn-Hilliard} is primarily due to tradition in the mathematical phase transitions community.}.  
 It may simply be viewed as a mathematical paradigm for {\it energy-driven pattern formation induced by competing short and long-range interactions}: Minimization of the first two terms  (short-range) leads to domains of pure phases of $u = \pm 1$ with minimal transition regions,   
whereas the third (long-range) terms induces oscillations between the phases according to the set volume fraction $m$. 
On a sufficiently large domain $\Omega$, the competition of the two leads to pattern formation on an intrinsic scale which depends entirely on $\gamma$. Throughout  this article, we always take the physical domain $\Omega$ of size\footnote{Since we do work on a finite domain (albeit sufficiently large), the choice of the exact domain size can still have an effect on the minimizing geometry, see Section \ref{domainsize}.}  much larger than this intrinsic scale ({\it cf.} Figure 
\ref{Fig:Intro}). 

A tool for our analysis is the $H^{-1}$ gradient flow ({\it cf.} \cite{CPW}) for $\bar u := u - m$ which  
takes the form 
\begin{equation}\label{pde}
\frac{\partial \bar u}{\partial t} \, = \, - \frac{1}{\gamma^2} \, {
\triangle^2} \, \,  \bar u\, \, + \,\,  
\triangle \left( { \bar u}^3 + 3 m { \bar u}^2 \,\,  - \, \, (1-3m^2){\bar u}\right) 
\, - \, { \bar u}. 
\end{equation}
It is important to note here that we compute the gradient with respect to  $H^{-1}$, a nonlocal metric. 
Hence  the presence of the  {\it nonlocal} term in the  functional  (\ref{OK}) simply gives rise to  a {\it local} perturbation\footnote{Note that the gradient of the nonlocal term with respect to $L^2$ would be $(-\triangle)^{-1} (\bar u)$, and  working in $H^{-1}$ has the effect of introducing an additional  $(-\triangle)$. } of the standard Cahn-Hilliard equation. 
However, it significantly changes its behavior, making it in some ways a hybrid of the Swift-Hohenberg equation \cite{SH} and the standard Cahn-Hilliard equation.

Whilst the functional  (\ref{OK}) is mathematically interesting on its own, 
we were drawn to it because of its  
direct connection with self-assembly of diblock copolymers: 
the functional is a rescaled version of a functional  introduced by Ohta and Kawasaki \cite{OK, NO}. 
Melts of diblock copolymer display a rich class of self-assembly nano-structures from lamellae, spheres, 
cylindrical tubes, to  double-gyroids  and other more complex structures (see for example, \cite{BF,Ketal}). Moreover, the usefulness of block copolymer melts is exactly this remarkable ability for self-assembly into particular geometries. For example, this property  can be exploited to create   
materials with   {\it designer}  mechanical, optical, and magnetic properties~\cite{BF}. 
Therefore from a theoretical point of view, one of the  
main challenges is to predict the phase geometry/morphology for a given set of material parameters, that is,  the creation of a phase diagram.   

The state of the art for predicting the phase diagram is via the self-consistent mean field theory (SCFT) \cite{MS, Fr}. 
While the simple  functional (\ref{OK}) can be connected with the SCFT via approximations ({\it cf.} \cite{CR}) with increasing validity close to the 
order-disorder transition, it is generally regarded as the basis of a qualitative theory, and one might question its usefulness with regard to predicting self-assembly structures for given material parameters. 
Preliminary numerical experiments  \cite{TN, CPW, TN2} indicate that 
all the phases  (including double gyroids and perforated lamellae), some of which had been predicted using the SCFT \cite{MS} and all of which have been
observed for polystyrene-isoprene \cite{Ketal}, can be simulated as minimizers of (\ref{OK})  starting from random initial conditions.  
This begs the question as to the extent to which a phase diagram via (\ref{OK}) can be compared with those of experimental observations  and SCFT calculations, at least close to the order-disorder transition. 

A thorough numerical phase diagram for (\ref{OK}) with $n=3$ (3D) is by no means an easy task.  In addition to numerical complications associated with the stiff PDE (\ref{pde}) and the necessary  small time steps and large spatial grids, the energy landscape of (\ref{OK}) is highly non-convex,  with multitudes of local minimizers and metastable states about which the gradient flow dynamics are very slow.
Many of the 
the simulations we presented in \cite{CPW}, were not simply the  final steady states for simulations of (\ref{pde})   with random initial conditions. Rather, one tended to get stuck in metastable states and some procedure for
{\it exiting} these metastable states in order to flow to {\it lower} energy states was crucial.  Such  procedures are often loosely dubbed {\it simulated annealing}. 
The present article may be viewed as a {\it test case} for the  3D phase diagram. Here, via a combination of  asymptotic analysis  and numerical experiments, we 
 consider the phase diagram in space dimension $n=2$  with respect to both {\it global} and {\it stable local} minimizers.
  The 2D situation is greatly simplified as the range of possible minimizers is both drastically reduced and, in fact, well accepted to be only certain basic 
structures: 
\begin{itemize}
\item 
disordered, i.e. minimizers of (\ref{OK}) are simply the uniform state $u \equiv m$; 
\item
lamellar, i.e. minimizers of (\ref{OK}) have a one dimensional structure; 
\item
spots, i.e. minimizers of (\ref{OK}) are a periodic array of (approximate) circles arranged on either  a hexagonal or rectangular lattice (in what follows, all spot solutions are hexagonally packed unless
otherwise indicated).   
\end{itemize}
Examples of the latter two cases and a metastable mixed state are presented in Figure \ref{Fig:Intro}.

We adopt periodic boundary conditions throughout this article and view the domain $\Omega$ as a flat torus.   
Our numerical and asymptotic results demonstrate  
explicit regions in the $\gamma$ vs $m$ plane for {\it stability} of the lamellar, circular and disordered states and 
explicit regions in  the $\gamma$ vs $m$ plane wherein the lamellar, hexagonally packed circular, and disordered states are {\it global minimizers}. At the surface, this may not to appear to be a particularly novel contribution. 
Indeed, 
there are many PDE models for pattern formation, some of which, like this one, are variational, and 
transitions between disorder, spots and lamellae as volume fraction changes and asymptotic-based descriptions are not surprisingly ubiquitous ({\it c.f.} \cite{PT, Ho}). 
For instance,  there is considerable interest in {\em localized} solutions to the Swift-Hohenberg equation and other similar equations in dimensions 1, 2 and 3 -- see for example \cite{SH, BK, SSM, Lloyd, PT, Ho, PF, Eletal} and the references therein. 
In \cite{PelWil} periodic solutions to the Swift-Hohenberg equation are analyzed in one dimension with both the detuning parameter $\alpha$ and the imposed fundamental period length as parameters. For our 
 problem,  the period length is determined by the minimization procedure, and not solely by a parameter. 
More importantly,  our focus here lies on the two-parameter, highly non-convex energy  landscape of the general functional  (\ref{OK}), not simply on steady state solutions to the partial differential equation (\ref{pde}). 
Specifically,   the  novelty of our approach resides in the:   
\begin{itemize} 
 \item[(i)] the demonstration of 
a hybrid numerical method  to address the complex metastability issues and 
 provide access to the ground state of (\ref{OK}); 
\item[(ii)]  
the (seemingly surprising) strong agreement (c.f. Section \ref{agree}) of our numerical results with the results of standard asymptotic analysis suggesting that  
for the purposes of describing the basic geometric morphology of the ground state of (\ref{OK}), the linear behavior of  (\ref{pde}) 
gives very accurate predictions in a (finite) neighborhood of the order-disorder transition
 \end{itemize} 
To this end, we are not aware of any systematic study for similar models with two parameters in two or more dimensions.  Finally, we note that both (i) and (ii) above  are particularly important with regard to a 3D study wherein the set of possible candidates for minimizers (global and local) is far more complex and in fact unknown.

 
\section{Asymptotic calculations}\label{asy}

In this Section we identify local and global minimizers of the energy functional (\ref{OK}) by determining the stationary solutions to (\ref{pde}) and evaluating their stability and energy. Asymptotic methods are used to describe the solutions in the limit $|\gamma-2| \ll 1$ and numerical 
experiments show they are valid over a much larger parameter region. We begin by
identifying the order-disorder transition (ODT) curve which separates regions of phase space where $u \equiv m$ is linearly stable from those regions where it is not.

\begin{figure}[t!]
\centerline{{\includegraphics[width =.5\textwidth]{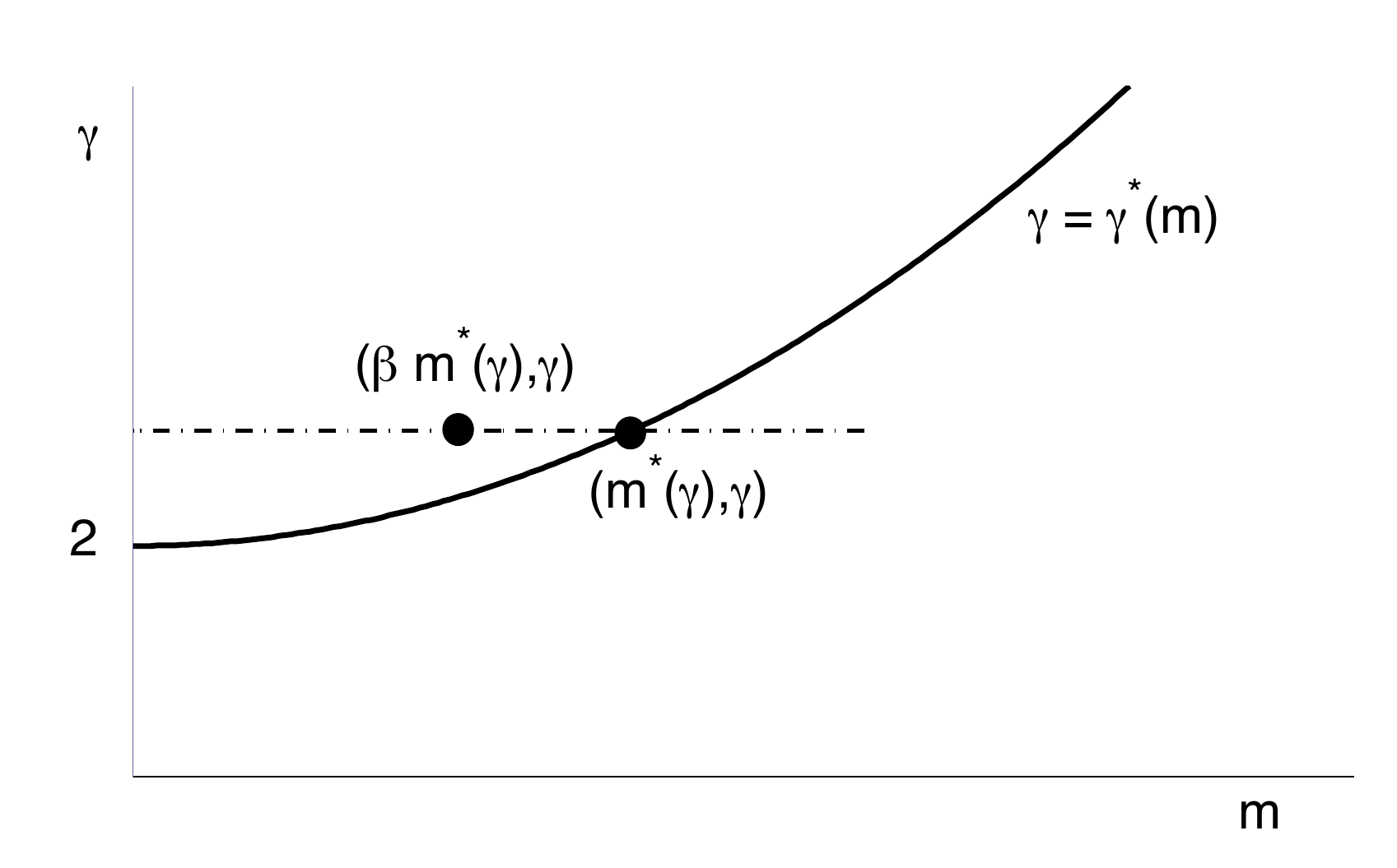}}}
\caption{Asymptotic coordinates.  We identify points in the plane via $(\beta m^*(\gamma),\gamma)$
where $m^*(\gamma) \ll 1$  as $\gamma \downarrow 2$.}
\label{fig:CoordsSketch}
\end{figure}

Consider the linearisation of (\ref{pde}) about $\bar{u} \equiv 0$
\begin{equation}\label{LDef}
v_t = \L v \equiv -\frac{1}{\gamma^2}\Delta^2 v - (1-3m^2)\Delta v - v
\end{equation}
where $v$ is periodic on $\left[-L_x/2,L_x/2\right]\times \left[-L_y/2,L_y/2\right]$. Taking the Fourier transform shows that
the eigenvalues of $\L$ are given by:
\[ \lambda = -\frac{1}{\gamma^2} |(k_x,k_y)|^4 + (1 - 3 m^2)|(k_x,k_y)|^2 - 1, \]
where $(k_x,k_y)$ represents the two-dimensional wavevector. The condition $\max(\rm{Re} (\lambda)) = 0$ determines the ODT curve as
\[ \gamma^* = \frac{2}{1-3m^2} \qquad {\rm or} \qquad m^* = \sqrt{\frac{\gamma-2}{3\gamma}}.\]
At fixed $\gamma$, $\bar{u} = 0$ is linearly stable for $m > m^*(\gamma)$ and unstable for
$m < m^*(\gamma)$.
On the ODT curve there is a finite dimensional center manifold on which we will construct solutions and determine their stability. 
For $ \beta >0$, we set 
$m = \beta m^*(\gamma)$ and  identify points in the plane via $(m,\gamma)$  
where $m^*(\gamma) \ll 1$  as $\gamma \downarrow 2$ (see Figure \ref{fig:CoordsSketch}). 
 Another of way of saying this is that  we consider the asymptotics as $( m, \gamma)$ approaches the point $(0, 2)$ along curves of the form $ \gamma \, = \, \frac{2}{1 - 3 \frac{m^2}{\beta^2}}.$ 
We introduce a regular asymptotic expansion for $\bar{u}$ for small $m^*$:
\begin{equation}
\bar{u}(x,y) \sim m^{*} \bar{u}_1(x,y) + (m^{*})^2 \bar{u}_2(x,y) + \ldots,
\end{equation}
where $\bar{u}_1 \ne 0$, since we are interested in the nonzero solutions of
(\ref{pde}). Consideration of the $\mathcal{O}(m^*)$ equation along with periodic boundary conditions leads to the following leading order representation of $\bar{u}$ for $m^* \ll 1$,
\begin{equation}\label{orderedstate}
\bar{u}(x,y,t) \sim m^*\big( a(t)\phi_1(x,y) + b(t)\phi_2(x,y) +  c(t)\phi_3(x,y)\big),
\end{equation}
where
\[\phi_1 = \sqrt{2}\cos{\left(\sqrt{2}x\right)}, \,
\phi_2 = \sqrt{2}\cos{\left(-\frac{1}{\sqrt{2}}x+\sqrt{\frac{3}{2}}y\right)},
\, \phi_3 =
\sqrt{2}\cos{\left(-\frac{1}{\sqrt{2}}x-\sqrt{\frac{3}{2}}y\right)}.\]
\begin{figure}[!t]
\center{
{\includegraphics[width=0.8\textwidth]{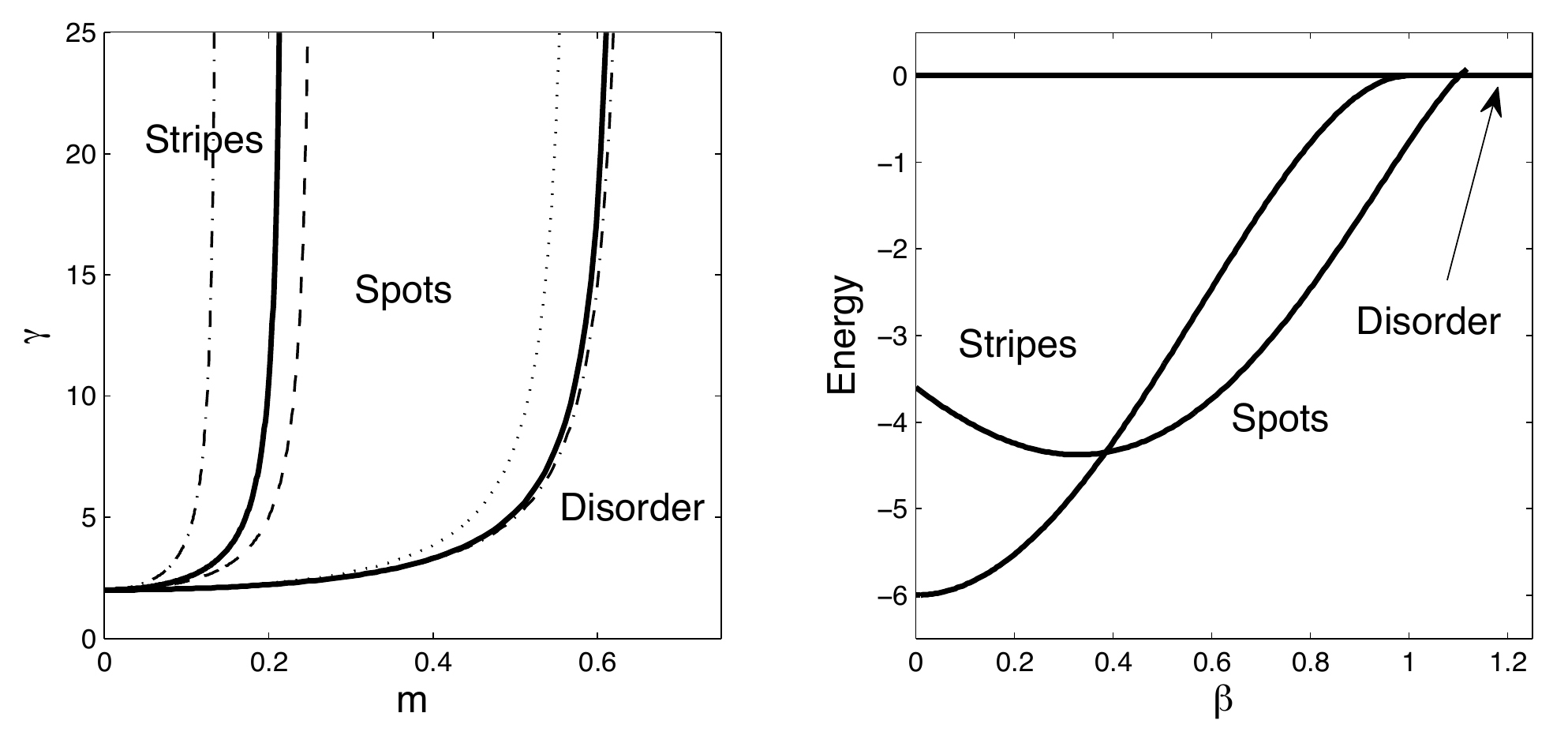}}}
\caption{Asymptotic results. Stability diagram of stationary solutions (left). Solid lines represent
global stability boundaries and dashed are linear stability boundaries. Lamellae are globally stable between $m=0$ and the first solid line and linearly stable until the dashed line. Spots are globally stable between the solid lines and linearly stable between the dash-dotted lines. 
The dotted line denotes the linear stability transition of the constant state. 
Energy (\ref{Lyap}) for steady solutions
of (\ref{amplitudesystem}) (right).}
\label{stab}
\end{figure}
\noindent
This expansion captures all the possible symmetric periodic configurations described above.
We proceed to determine the amplitude dynamics on
the center manifold of (\ref{LDef}) $X^c=\mathrm{span}\{\phi_1,\phi_2,\phi_3\}$ by projecting
the full PDE  (\ref{pde}) onto $X^c$. Thus, we consider
\[
\left \langle \bar{u}_t+\frac{1}{\gamma^2} \, {\triangle}^2 \bar{u} - 
\triangle \left( \bar{u}^3 + 3 \beta m \bar{u}^2\right) + (1-3\beta^2m^2){\triangle
  \bar{u}} + \bar{u}, \phi_i\right \rangle=0,
\]
for $i=1,2,3$ on the domain
$\Omega$, appropriately defined by
$L_x = \frac{n\pi}{\sqrt{2}}$ and $L_y = \frac{2}{\sqrt{3}}L_x$. Computing the inner
products  in $L^2$ and expanding in powers of $m^{*}$ with the aid of Maple, we find
the following amplitude ODE system
\begin{eqnarray}\label{amplitudesystem}
\dot{a} & = & 6(1-\beta^2)a - 6\sqrt{2}\beta bc - 6(b^2+c^2)a - 3a^3 \nonumber \\
\dot{b} & = & 6(1-\beta^2)b - 6\sqrt{2}\beta ac - 6(a^2+c^2)b - 3b^3  \\  
\dot{c} & = & 6(1-\beta^2)c - 6\sqrt{2}\beta ab - 6(a^2+b^2)c - 3c^3 \nonumber 
\end{eqnarray}
where time has been rescaled with $(m^*)^2$.

The original equation is a gradient system, and thus, the
reduced system is one as well. The corresponding Lyapunov function
is given by
\begin{eqnarray}
V(a,b,c) & = &  -3(1-\beta^2)(a^2+b^2+c^2) \,\, + \, \,  6\sqrt{2}\beta abc \nonumber \\
& & \quad + \,\, 3(a^2b^2+b^2c^2+a^2c^2)\, \, +\,\,  
\frac{3}{4}(a^4+b^4+c^4).
\label{Lyap}
\end{eqnarray}
 This function is in fact the projection of the energy functional (\ref{OK}) onto the center manifold $X^c$.
We consider the structure and stability of the stationary
solutions of the amplitude ODE system by analyzing its Lyapunov
function. The stationary states satisfy the following system
\begin{equation}\label{Vsystem}
V_a(a,b,c)=0, \, \, V_b(a,b,c)=0, \, \, \mathrm{and} \, \,
V_c(a,b,c)=0,
\end{equation}
and they are linearly stable when all eigenvalues of the
Hessian matrix, $H(V(a,b,c))$, are positive.   
We identify five fixed points of the system (\ref{Vsystem}):
\begin{enumerate}
\item {$a=0, \, b=0, \, c=0$  $\rightarrow$ \it{disorder};}
\item {$a=\pm \sqrt{2(1-\beta^2)},\,b=0, \, c=0$  $\rightarrow$
    \it{lamellae};}
\item {$a=b=\pm \sqrt{\frac{2(1-\beta^2)}{3}}, \, c=0$
    $\rightarrow$ \it{Triangular spots};}
\item {$a= \bar{a},\, b=c=\pm \bar{a}$, where  $\bar{a}=-\frac{\sqrt{2}}{5}\left[\beta \pm \sqrt{5-4\beta^2} \right] \rightarrow$ \it{hexagonally packed circular};}
\item {$a=b=\pm \bar{a} , \, c=-\bar{c}$, or $a=-b=\pm \bar{a}, \, c=\bar{c}$, where
    $\bar{a}=\sqrt{\frac{2(1-5\beta^2)}{3}}$ and $\bar{c}=2\sqrt{2}\beta$
    $\rightarrow$ $|a|=b\ne c$ {\it case}.}
\end{enumerate} 
While other similar systems can generate other patterns ({\it cf.} \cite{Ho}), none of these appear in this system.
Evaluating the eigenvalues of $H(V(a,b,c))$ at the five fixed
points, we determine the following three regions of linear stability:
\begin{enumerate}
	\item {$0 \le \beta < \frac{1}{\sqrt{5}}$: {\it lamellae};}
	\item {$\frac{1}{\sqrt{17}} < \beta < \frac{\sqrt{5}}{2}$: {\it hexagonally
	      packed circular};}
	\item {$\beta > 1$: {\it disorder}.}
\end{enumerate} 
Evaluation of the Lyapunov function at the three linearly stable
steady states is then used to determine the following three regions of
global stability: 
\begin{enumerate}
\item {$0 \le \beta < \frac{1}{29}\sqrt{551-174\sqrt{6}}$: {\it lamellae};}
\item {$\frac{1}{29}\sqrt{551-174\sqrt{6}} < \beta < 3 \sqrt{\frac{5}{37}}$: {\it hexagonally
      packed circular};}
\item {$\beta > 3 \sqrt{\frac{5}{37}}$: {\it disorder}.}
\end{enumerate}

Note that there are regions where more than one solution type is linearly stable.
Figure \ref{stab} shows the asymptotic linear stability and global stability
diagrams.

\section{Numerical Simulations}
\begin{figure}[th!]
	\centerline{
{\includegraphics[width=.85\textwidth]{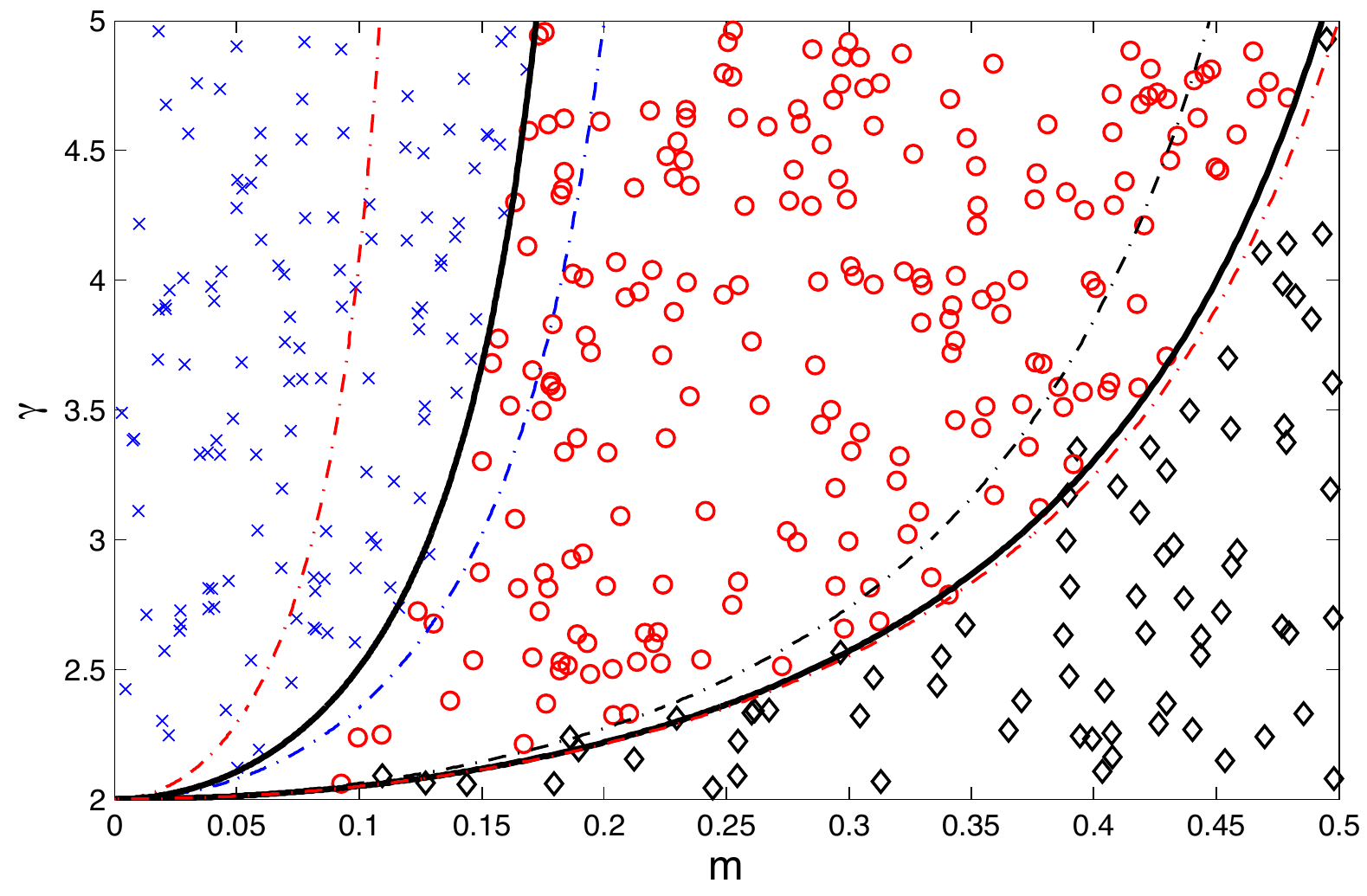}}}
\centerline{
{\includegraphics[width=.85\textwidth]{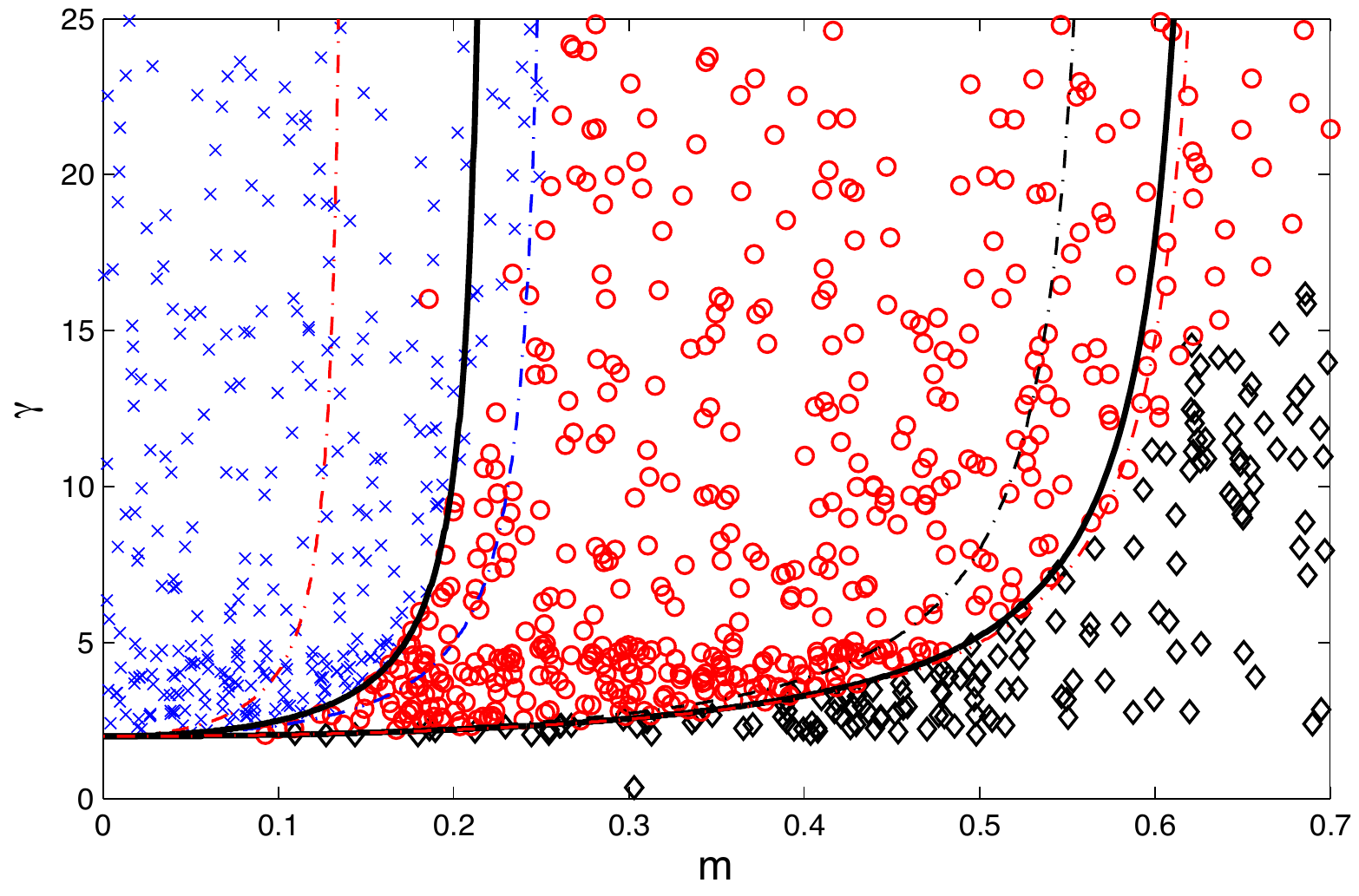}}}
\caption{
(These figures need to be viewed in colour) Numerically computed phase diagram. (Bottom) Complete diagram. (Top) Detail for  $\gamma$ close to 2. Blue crosses: Lamellae. Red circles: Hex packed spots. Black diamonds: disorder. The red dashed-dotted lines mark the linear stability boundary of spots, the blue dashed-dotted line marks the linear stability boundary of lamellae, the black dashed-dotted line marks the linear
stability boundary of the disordered sate,  and the solid black lines mark the global stability regions of lamellae and spots respectively. 
}
\label{Fig:GammaBif}
\end{figure}

In this section, we present the results of our numerical experiments and compare them with the asymptotic calculations of the previous section in regimes where 
they may present some validity. Our  numerical method will be  described in some detail in Section \ref{details}. It is a hybrid method which not only integrates the PDE 
(\ref{pde})  but also involves an interplay with the energy (\ref{OK}) through certain methods of {\it simulated annealing}. For most values of the parameters $(m, \gamma)$, we believe that this method results assesses  the ground state of (\ref{OK}),  at least from the point of the it's inherent structure (geometry and  symmetry). Of course we have no proof of this statement, and moreover, rigorous results pertaining to the ground state are difficult to obtain even in 2D (see \cite{CO, Mu2, ST, Sp} for partial results).  

We sampled the parameter plane by taking 
a sequence of randomly chosen points $(m_i, \gamma_i) \in [0, 1]\times[2,25]$ and implemented our method for 
 each such $(m_i, \gamma_i)$ 
 with $\bar{u}(x_i,y_j, t = 0) \in (-1,1)$ randomly chosen from a uniform
distribution. Once a sequence of more than 500 runs was complete, an edge detection algorithm was used to place additional points near the interfaces between regions where $\bar{u} = 0$, hex spots
 and lamellae are stable. For all runs  we  computed the energy  (\ref{OK}) at each time step, making sure that the final presented state had the least energy over the course of the entire run.

\begin{figure}[h!]
	\centerline{
{\includegraphics[width=.5\textwidth]{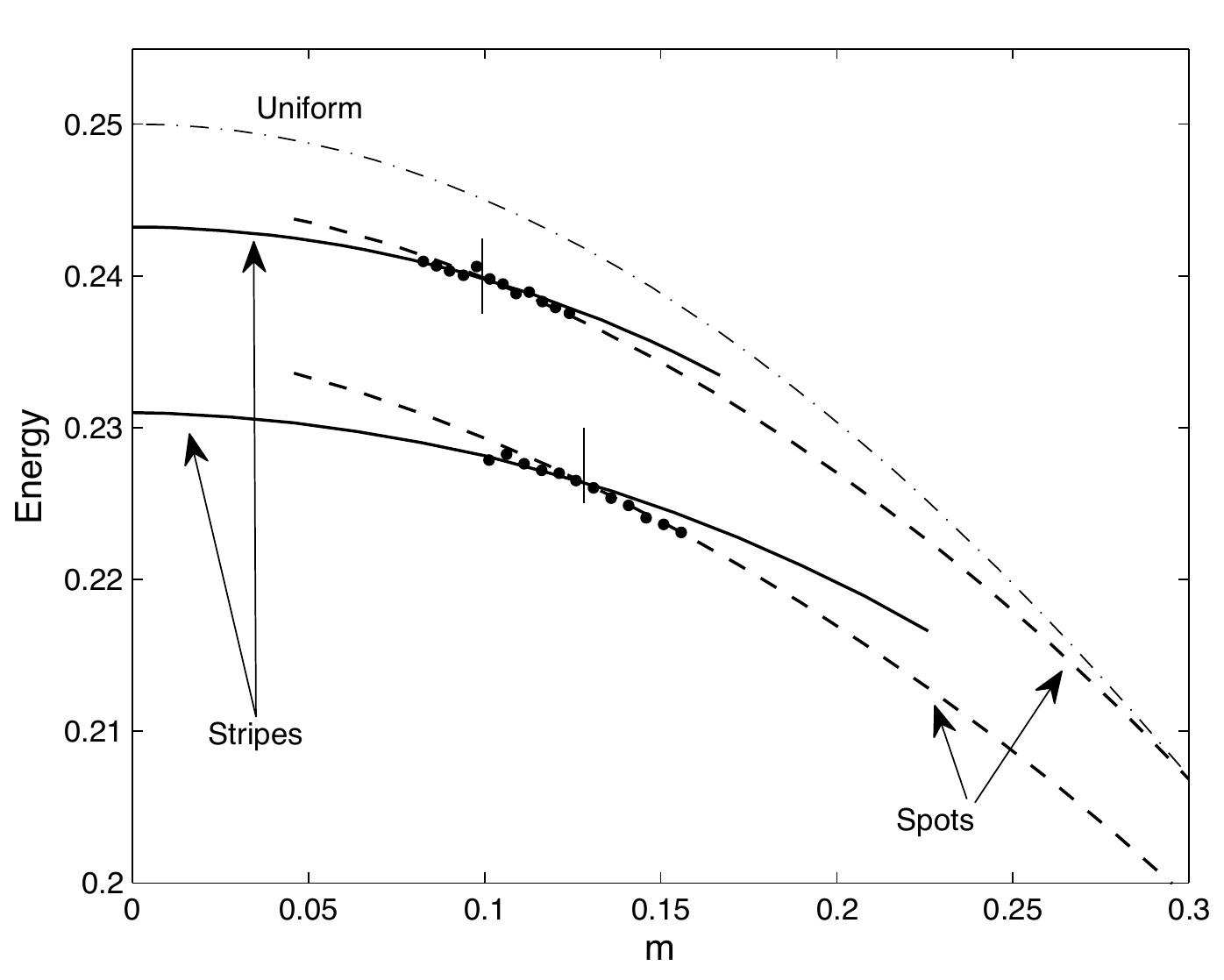}}
{\includegraphics[width=.5\textwidth]{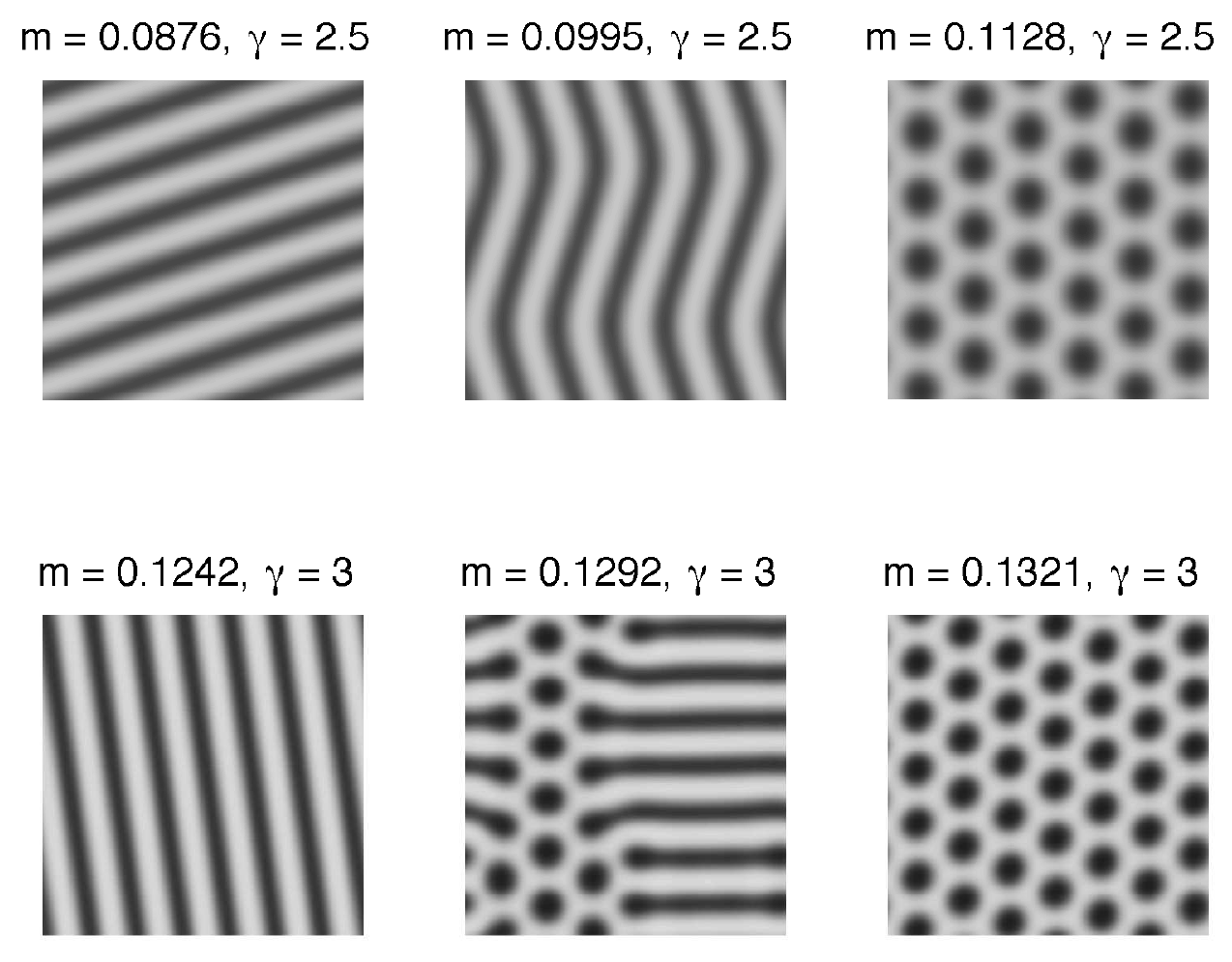}}
\vspace{.5cm}}
\centerline{
{\includegraphics[width=.5\textwidth]{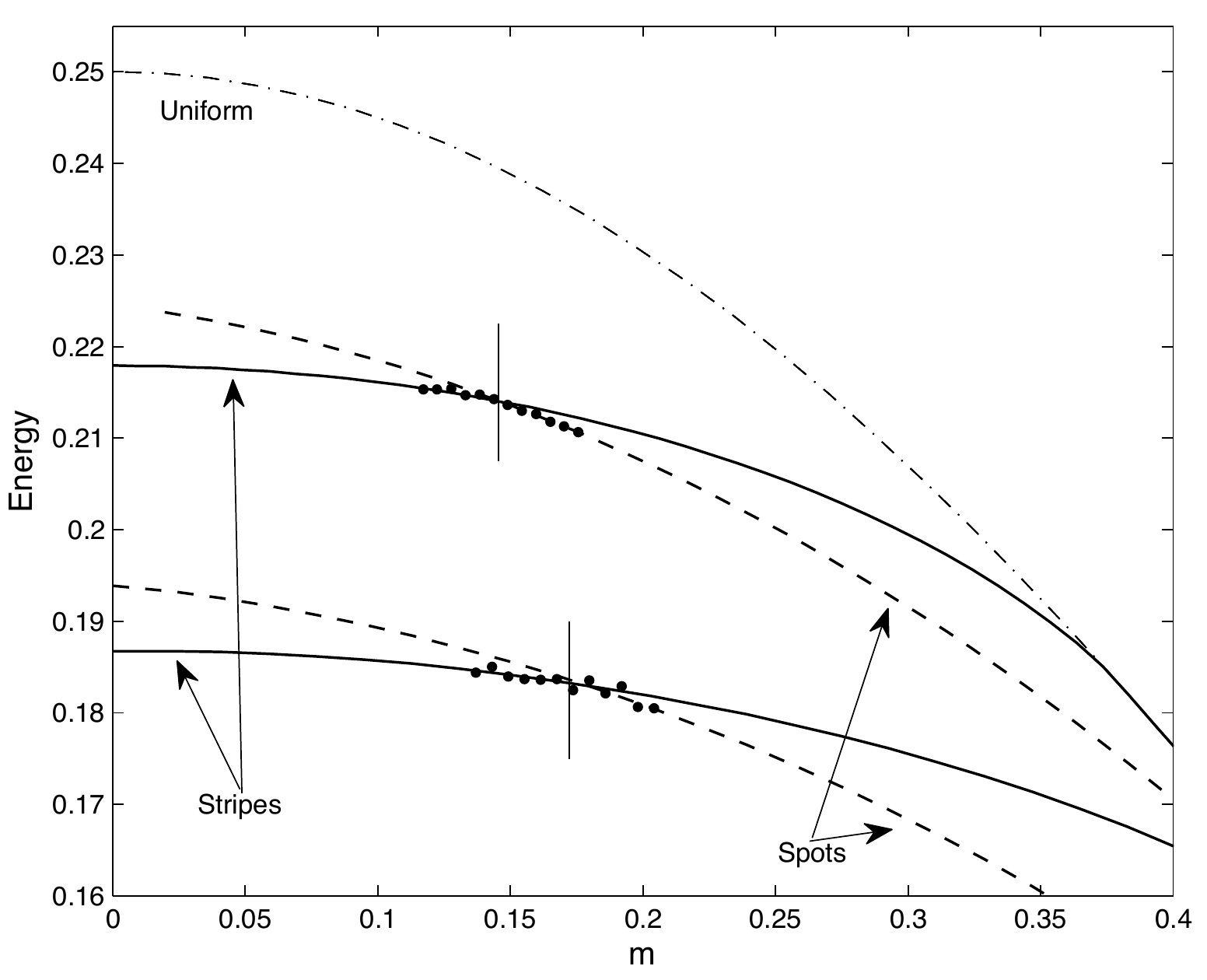}}
{\includegraphics[width=.5\textwidth]{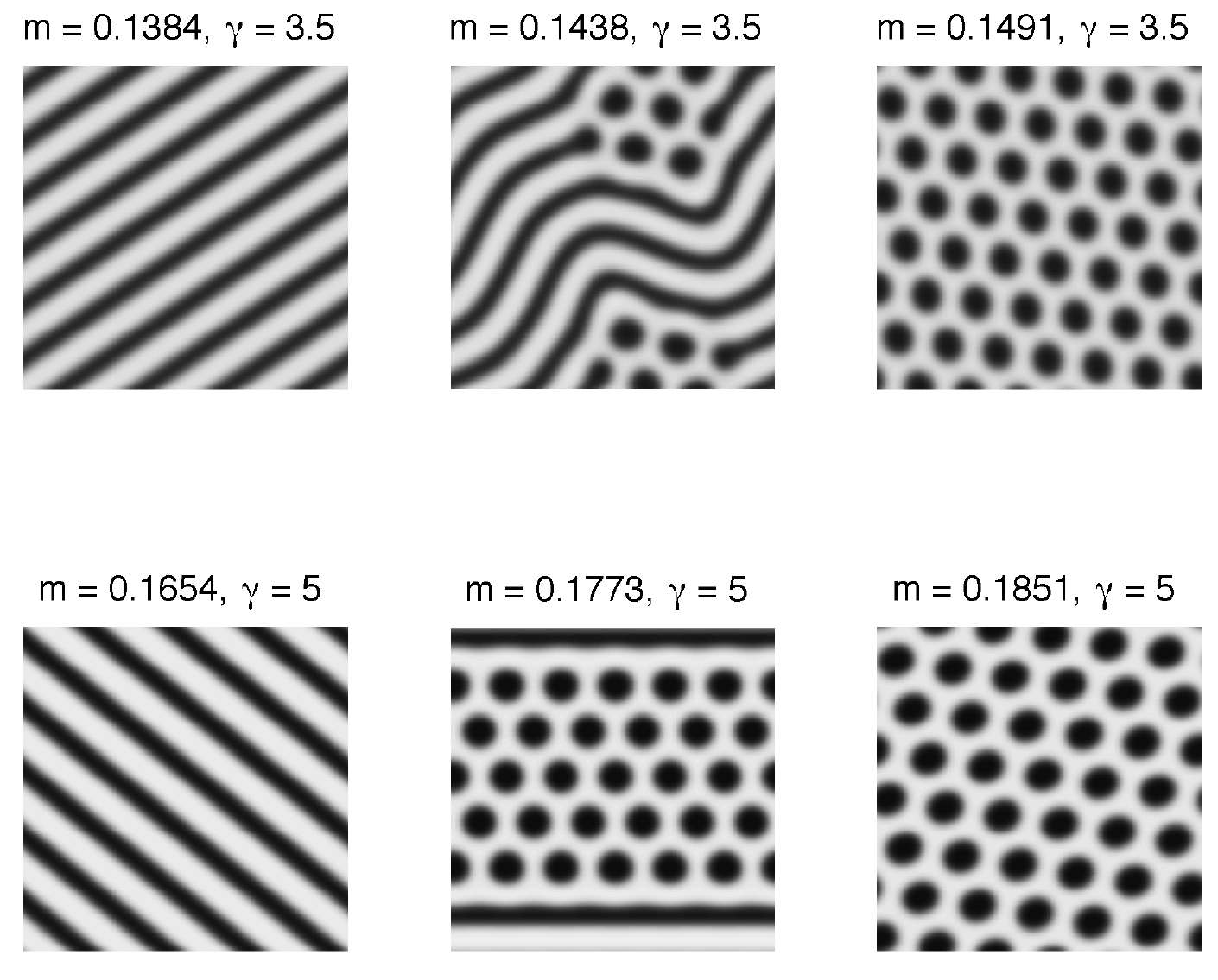}}
\vspace{.5cm}}
\centerline{
{\includegraphics[width=.5\textwidth]{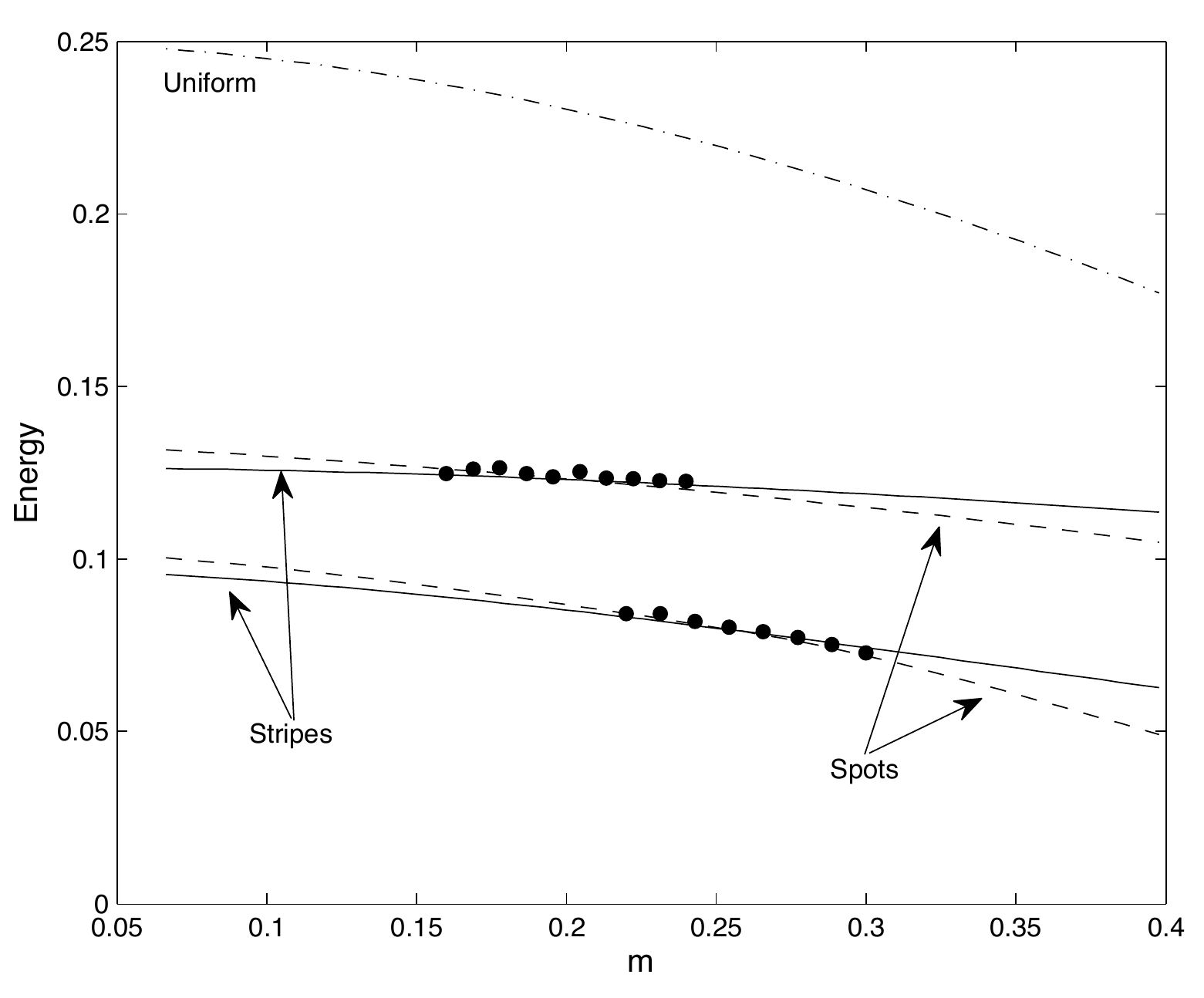}}
{\includegraphics[width=.5\textwidth]{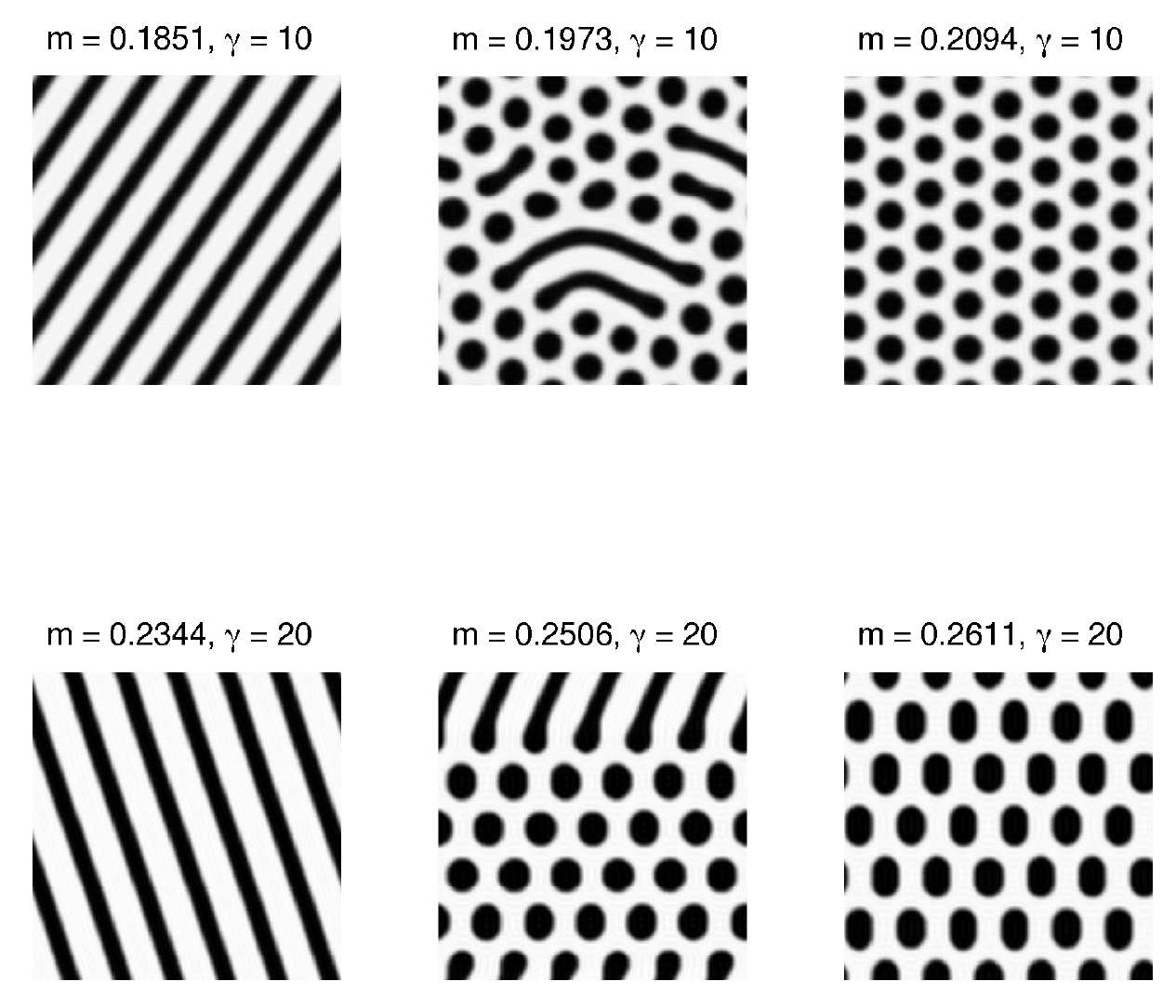}}
}
\caption{Comparison of energies along solution branches for $\gamma = 2.5, 3, 3.5, 5, 10$ and 20. 
Note that the
continuation routine failed for some values of $\gamma$ on the stripe branch for $m$ sufficiently large and on the spot branch for
$m$ sufficiently small but these always occurred well beyond the value of $m$ at which the branches
exchange global stability. In each energy diagram, the lower value of $\gamma$ has the higher energy.  The vertical bar is the asymptotic value $m^*(\gamma)$ where the branches exchange stability. On the right,
there are three plots for each run from random initial data: $m < m^*, m \simeq m^*, m > m^*$.  Here we can see that there is no convergence sufficiently close to $m^*$ but the expected patterns appear away from it.
Notice that the there is greater separation in the colors and hence the solution extremes as $\gamma$ increases. }\label{Fig:ContComp}
\end{figure}
\afterpage{\clearpage}
\subsection{Results of numerical simulations}

Our numerical experiments confirm the asymptotic description of the stability curves and the solution structure
in the limit $\gamma \downarrow 2$, $m^*(\gamma) \downarrow 0$. In Figure \ref{Fig:GammaBif} we present the numerically computed 
bifurcation diagram with the asymptotic stability curves overlaid. Here $\times, \circ$ and $\diamond$ indicate
stripes, hex packed spots and the disordered state respectively. The solid lines mark the global stability curves and the dashed lines outline the linear stability regimes. Note that almost no runs converged to a state outside of its asymptotic region of global stability and no runs converged to a state outside of its region of linear stability. Figure \ref{Fig:GammaBif} (bottom) displays a detailed examination of the bifurcation diagram for $2 \le \gamma \le 5$. This is within the expected range of asymptotic
validity and shows very good agreement.

Figure \ref{Fig:GammaBif} (left) has $2 \le \gamma \le 25$ and the 
asymptotic estimates are now only in qualitative agreement as $\gamma$ increases. Limited additional runs up to $\gamma = 100$ were
also performed and still show qualitative agreement with asymptotics in that no new phases were seen and the
ordering of stripes, spots and the homogeneous state occur for increasing $m$.
Unfortunately, the stability region of the PDE time-stepper quickly decreases with increasing $\gamma$ and the time-scale of PDE evolution slows down making the runs with very large $\gamma$ increase in cost faster than $\gamma^2$ and hence very computationally expensive.  For increasing $\gamma$ we find that the stripe/spot
transition is much better approximated by the asymptotics than the spot/homogenous transition.  That is, for
increasing $\gamma$ there is an ever wider region where spots are globally stable but the homogeneous 
state is linearly stable.  

Continuation in $m$ of spots and stripes was also performed
for $\gamma =$  2.001, 2.01, 2.1, 2.25, 2.5, 3, 3.5, 5, 10 and 20.  Comparing the energies of the different phases allows us to identify the global
minimizer directly.  Figure \ref{Fig:ContComp} shows that the agreement with the direct PDE computations is very good. However, computationally continuation is much slower and is not
completely reliable.  Continuation keeps the geometry fixed as $m$ varies not necessarily
guaranteeing a global minimizer over the whole of the branch.  Typically, we computed three 
branches of solutions and kept the ones with lowest energy in the region of the transition.

These results clearly demonstrate that the asymptotic analysis presented in Section 2 accurately describes both the phase diagram and solution structure of minimizers beyond the point
$(m,\gamma) = (0,2)$ and up to perhaps as far as $\gamma = 5$.  For $\gamma > 10$, the agreement is
much weaker between the asymptotics and numerics.

\section{Numerical methods}\label{details} 
The details of many  aspects of our numerical implementation are straightforward and well understood. We use a pseudo-spectral method in space, because of the large scale periodicity, and two different timestepping schemes for evolution. For early times, we use
exponential time-differencing (ETD) as the dynamics are quick and we need small timesteps to resolve them.  ETD provides a cheap per step highly accurate method with well known
stability and accuracy properties for stiff diagonalizable PDEs.  For later times we switch to an 
iterative linearly implicit gradient stable algorithm \cite{Xu}.  This method is more expensive per
step and only first order accurate but it allows arbitrarily large time steps and guarantees that the
energy will not increase.

For all computations we used $N=256$ spatial modes in both dimensions and $\Delta t = \frac{.1}{1+\gamma^{3/2}}$. 
The initial choice  of domain size was somewhat arbitrary ({\it cf.} (\ref{domain-size})). As we have mentioned the domain size used is sufficiently large enough with respect to the intrinsic period scale which is determined by $\gamma$. However, even with such a choice, the exact size of the domain can influence both 
the minimizing patterns  and certainly the gradient flow route towards them. Thus we have accounted for this in our algorithm via a {\it variation of domain size}. 
This is discussed in Section \ref{domainsize}.   
Typically we computed for $0 \le t \le 100 = t_F$ but took $t_F$ as large as 2500 in some cases. All computations were done on a laptop and implemented in Matlab. Unless otherwise indicated, the norm of the residual ($||u_t||_2$) is less than
$10^{-8}$ in all figures.

\subsection{Metastability}\label{sec-meta}
\begin{figure}[!t]
\centerline{{\includegraphics[width=0.4\textwidth]{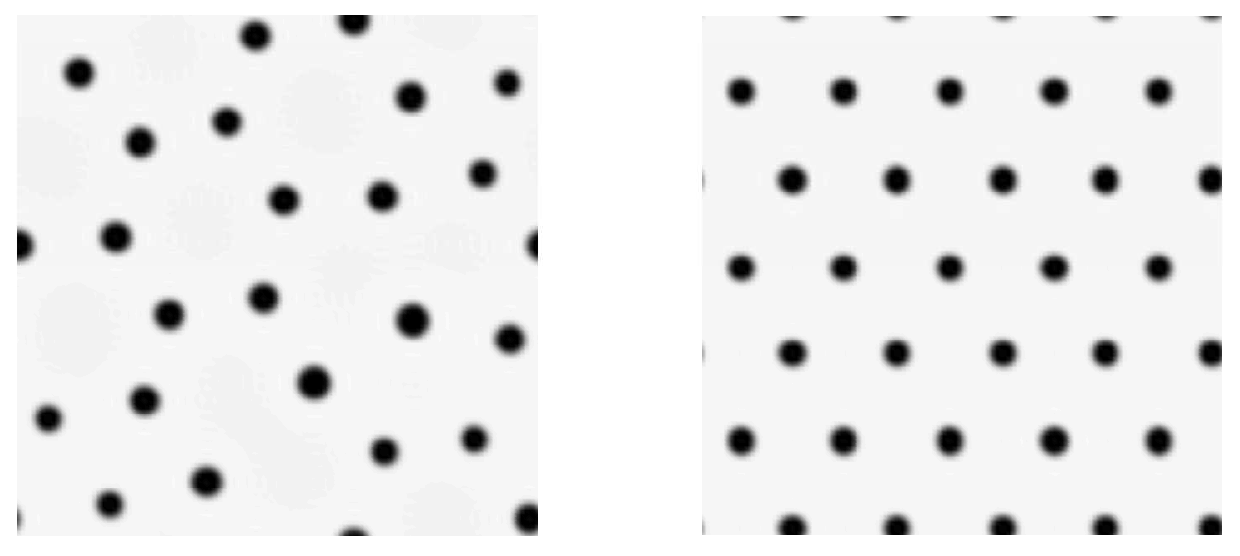}}}
\caption{Metastable state with $\gamma = 20, m = 0.8$. Left: $t = 100$. Right: $t = 10000$. The difference in energies between these profiles is less than $1e-7$ but clearly the left state is not an energy minimizer as the packing of spots is irregular. The dynamics is driven by weak interaction between the spots.}
\label{Fig:Meta1}
\end{figure}

\begin{figure}[!t]
\centerline{{\includegraphics[width=0.8\textwidth]{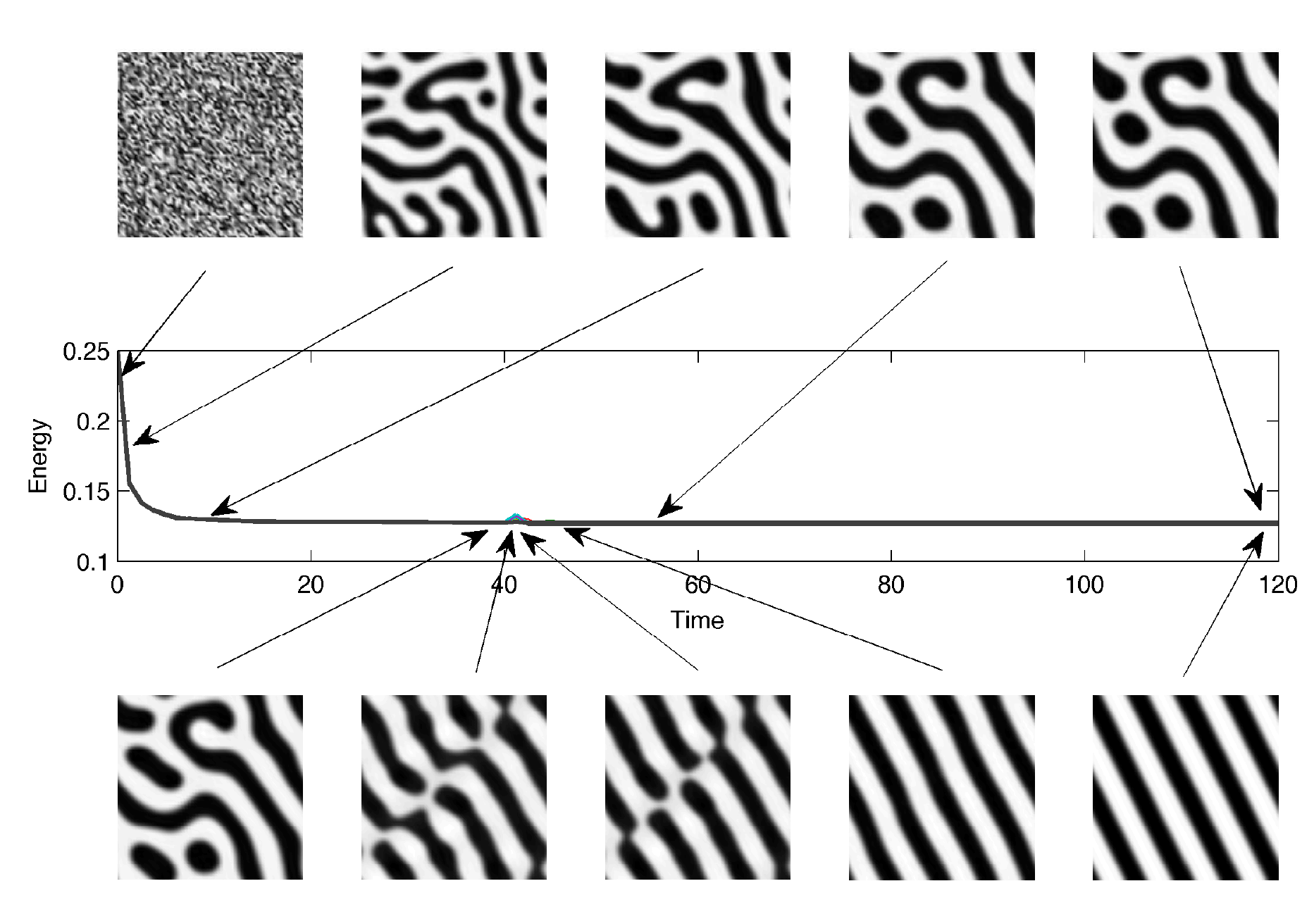}}}
\caption{Comparison of two runs with the same initial conditions. The top simply integrates the PDE while the bottom implements the spectral weighting algorithm described in the text. Both runs are identical for $0\le t < 40$.  At $t=40$ the spectral weighting is turned on and acts very quickly to allow lamellae to develop. Note that at this scale the difference in the energies is not noticeable (see Figure \ref{Fig:Energy} below). In both runs $\gamma = 10$ and $m = 0.1$.}
\label{Fig:Meta}
\end{figure}
\begin{figure}[!t]
\centerline{{\includegraphics[width=0.8\textwidth]{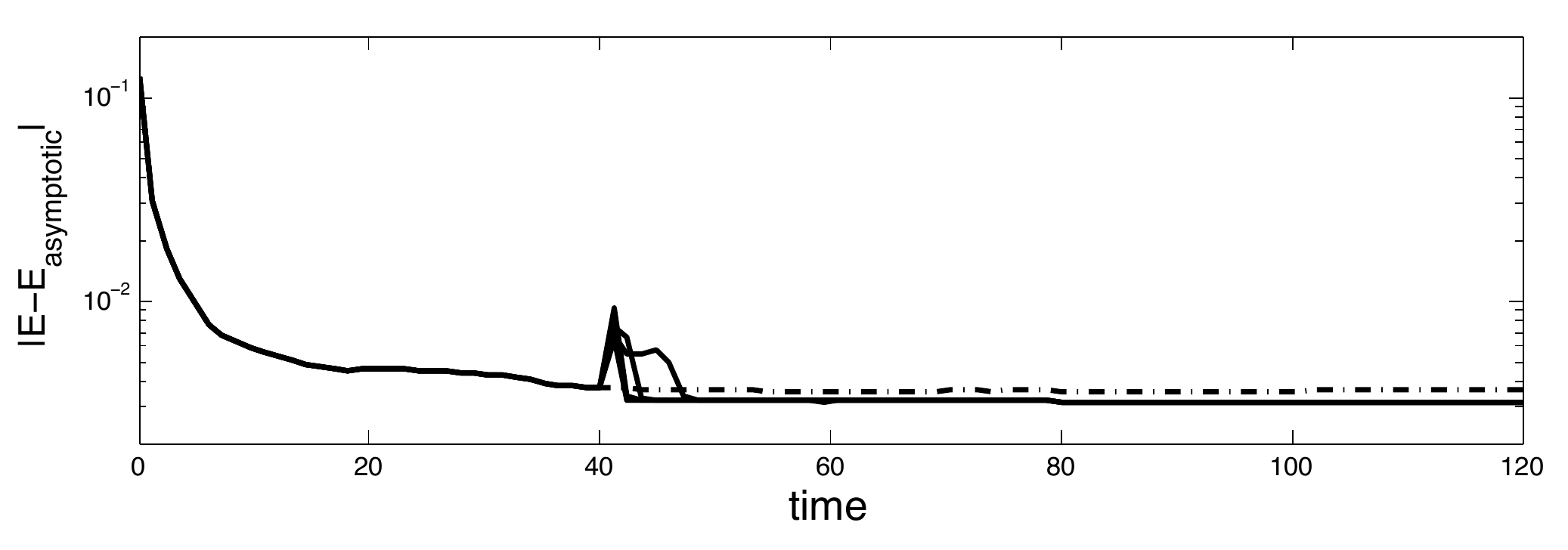}}}
\caption{Detail of energy over time for Figure \ref{Fig:Meta}.  The dashed line has no spectral weighting. The solid lines all have spectral weighting for $40 < t < 80$. The damping parameter $\rho$ takes the values $.05, .1, .15, .2, .25 $ and $.3$. Over this range this is no discernible 
difference in the final states with weighting. For all runs $\gamma = 10$ and $m=0.1$.}
\label{Fig:Energy}
\end{figure}

Figure \ref{Fig:Meta} shows snapshots of a typical run as well as the energy decay over time.  Notice that the ``final" profile is a labyrinthine-like pattern rather than pure lamellae or spots but that the energy decay is, relatively, small at that time. This phenomenon is pervasive for this problem, as the functional ({\ref{OK}) has many local minimizers and metastable states near which the
dynamics are very slow. Numerically, one cannot distinguish a metastable state from a stable
one, since they are both identified as solutions for which the relative
change in $\bar{u}$ or $E$ between timesteps is smaller than some tolerance level. Furthermore,
irregular long lasting states are common in diblock copolymer
experiments \cite{Fr}, so, without additional analysis, it is unclear whether they
are steady states, rather than just persistent intermediate
profiles.

Analytical considerations can provide insight into the classification of
observed profiles into metastable and stable steady states.  The analytical results on the energy functional (\ref{OK}) 
are {\it suggestive} that its global minimizers are {\it nearly} periodic profiles with level surfaces with {\it nearly} constant mean curvature  
  \cite{ACO, Sp, Mu, CO, Mu2, ST, CP2}. The 
irregular states that we observe do not match this description, and
thus we term them as metastable states. Long-lived transient behavior due to the  metastability of the Cahn-Hilliard
equation is well-known (see for example, \cite{Ward}), and the metastability is also observed
in  SCFT numerical calculations \cite{Fr}. Given that our PDE differs from the
Cahn-Hilliard equation only in the extra nonlocal term, and can be connected to the SCFT theory 
(albeit, via  further approximations), there is little reason to
believe that it can escape the issue of metastability in general.  Hence, we need
to modify the gradient descent to better approximate the true global
minimizers.

\subsection{Selecting and damping modes}

Techniques for dealing with metastability and highly non-convex energy landscapes often belong to the broad
class of statistical methods, called simulated annealing.  They  were created 
to navigate through a complex energy landscape
in search of a global minimizer. A very simple form of simulated
annealing can be achieved by adding unbiased noise to the evolved metastable
state. This may force the solution out of the local minimizer that it
is stuck in and make it continue its evolution through the energy
landscape. Unfortunately, this approach does not provide a guaranteed way of
addressing metastability as: too much noise leads to the divergence of
the solution and even when the solution remains bounded, there is no
way of ensuring that it will not revisit the local minimizers that it
was stuck in before. Also, the added noise is very quickly damped out due to the
fourth-order derivative term in this problem.  

A different approach to the removal of defects is
provided by the technique of spectral filtering. The essence of this
method lies in the removal of insignificant spectral components from
an evolved state. That is, we evolve the solution
from random initial conditions
until a structure is formed, compute its Fourier coefficients, and keep only the modes which correspond to
the coefficients above a certain threshold. The evolution is then
continued and the process of spectral filtering repeated. This
approach was suggested in \cite{Fr}, and, combined with adding noise, applied to our PDE in a 
modified form.

Functional (\ref{OK}) can be thought of as a 
length selection mechanism and, in fact, if we plot the energy concentration in Fourier space over time we see that this happens quite quickly. Figure \ref{Fig:EnergyConc} (left) shows $\hat{v}(k^*,t)$ such that $\max_k \hat{v}(k,t) = \hat{v}(k^*,t)$. In Figure \ref{Fig:EnergyConc} (center)  we see that the energy is concentrated in a narrow band about $k^*$ in Fourier space.  Motivated by this
observation we evolve the original PDE until past the point that the dominant length-scale has
emerged. Denoting this length scale $k^*$ we damp the Fourier coefficients as follows: 
$\hat{v}(k) \to w(k;k^*) \hat{v}(k)$ where 
\begin{eqnarray*}
w(k;k^*) &=&     (1-\rho) \, + \,  \\
&& \hspace{-1cm} \rho \left( \exp\left({-5(1-|k|/k^*)^2}\right)
\,  + \, \exp\left({-5(2-|k|/k^*)^2}
\right) \, + \, 
\exp\left({-5(3-|k|/k^*)^2} \right) \right)
\end{eqnarray*} 
This keeps information at all wavelengths but focusses the dynamics at the key length-scale and
its higher harmonics. Experimentation with the parameter $\rho$ indicates that there is little difference in the outcome with $0.05 \le \rho \le 0.3$ as indicated in Figure \ref{Fig:Energy}. With $\rho$ too small there 
is no effect and the wrong pattern may emerge if $\rho$ is too large or the standard gradient flow
is not run long enough.

This approach allows us to first identify the energy minimizing length scale, the dynamics to
focus on these wave lengths, ensure the local stability of profile which emerged and then finally to smooth out the effects of added noise. Figure \ref{Fig:Meta} shows two runs from the same initial
conditions. The top has no spectral damping and ends up stuck in a metastable state whereas the
bottom leads to the global minimizer. The energy of both runs is shown in the middle Figure \ref{Fig:Meta} where the difference between the two runs is almost indistinguishable.  A detailed view of the energy is presented in Figure \ref{Fig:Energy}. Lastly, Figure \ref{Fig:EnergyConc} (right) shows the coefficient distributions of the final time profiles in Figure \ref{Fig:Meta} with $+$ for global minimizing state and $\cdot$ for the metastable one. Notice that the spectral weighting did not shift $k^*$ but rather allowed a simpler pattern to emerge.

\begin{figure}[!t]
\centerline{{\includegraphics[width=\textwidth]{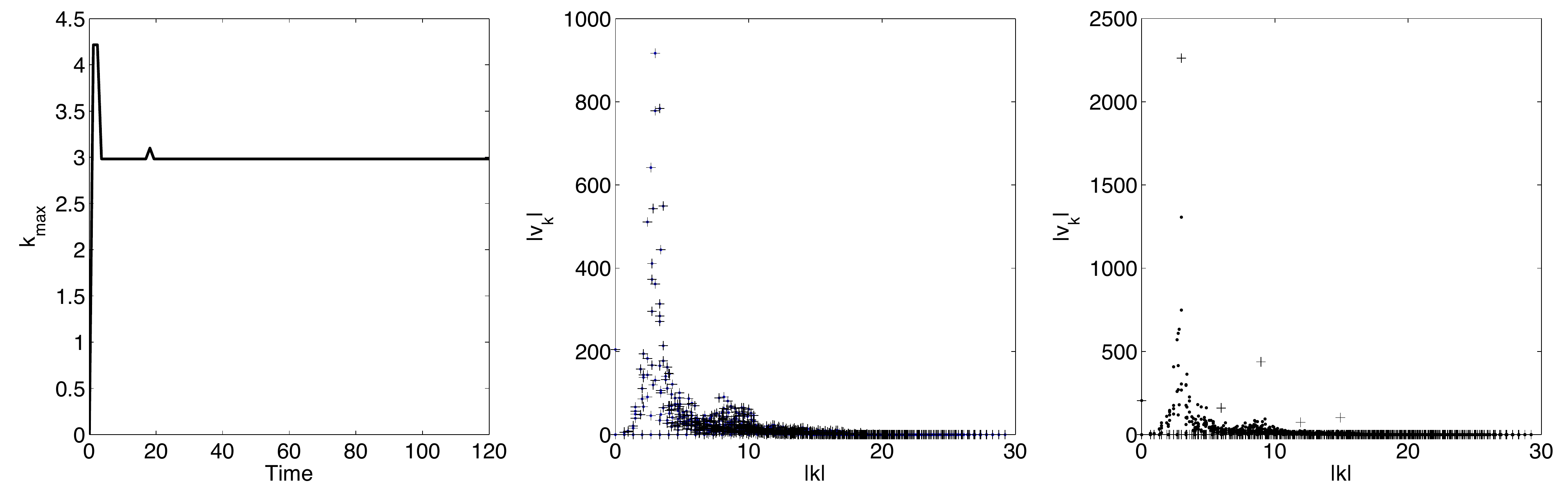}}}
\caption{(left) $k^*$ over time. Here we see that the evolution quickly converges to a dominate length scale and that the spectral weighting does not change it. (middle) before the application of the spectral weight most energy is concentrated near the dominant mode $k^*$. (right) The effect of the spectral weight is to allow the energy to concentrate even more so on $k^*$ and its integer multiples.}
\label{Fig:EnergyConc}
\end{figure}

\vspace{.5cm}
\noindent
{\bf Remark}
{\em We have also included four movie files that clearly demonstrate the use of this approach.  We compare two cases $(m,\gamma) = (.1,10)$ and $(m,\gamma) = (.25,10)$ showing the effect of taking $\rho = 0$ or $\rho = .1$. For $m = .1$, the unmodified run leads to a metastable mixed state including both spots and stripes whereas modified run leads to pure stripes with lower energy. For $m = .25$ the unmodified run leads to a non-optimal collection of spots of differing sizes distributed seemingly randomly. Here, the modified run leads to uniform spots packed on a hexagonal grid. All four runs were started from the same random initial data.}

\subsection{Choice of time-stepping method}
When deciding on the most appropriate time discretization scheme for numerical simulations of the PDE, we considered exponential time-differencing methods of different order and construction as well as
two gradient stable schemes based on the work of Eyre \cite{Ey}. We settled on the ETDRK4 scheme \cite{Tref} as it
has the highest order of accuracy and its stability region is larger than that of the ETD2 scheme. However, the time step is stability limited so to compute for long times when the dynamics are slow we need something else. For large timesteps we found that a fully implicit gradient stable method \cite{Ey,Xu} was
the most robust in general and the fastest when we could take large timesteps. Combining the two methods
led to an algorithm at least 3 times faster than either one alone.

The gradient stable scheme is stable for all timesteps and guarantees decay of the energy but requires 
a nonlinear solve at each timestep.  The authors in \cite{Xu} developed an iterative scheme to do this
but the number of iterations required greatly depends on the magnitude of the solution change over each
timestep.  This means it is not feasible to use it to take very large steps at the beginning. To
further speed the computation we developed an adaptive timestep routine based on the number of required iterations.


We also considered numerical implementation of the $L^2$ volume constrained gradient flow of the energy functional (\ref{OK}). This consideration was motivated by the fact that this evolution system is less stiff than the PDE (\ref{pde}), as it contains a laplacian, rather than a bilaplacian term. However,  we found that it did not result in a more efficient computation of the steady states as it required far more computational steps. This could be because the lower order derivatives do not damp out high frequency modes as quickly.

Continuation in $m$ was done with a modified Newton's method allowing for linearly implicit iterations
without computing the $256^2\times256^2$ Jacobian matrix.

\subsection{Domain size selection}\label{domainsize}
We do not know the natural period length of arbitrary solutions and we do not want to enforce a 
solution type due to our choice of box size thus we adaptively choose the domain size. In the limit
$\gamma \downarrow 2$ the asymptotic form of the solution has $|k| = \sqrt{2}$ so we take
$L_0 = 12 \pi$ as a preliminary box size to fit a reasonable number of periods in the box. As $\gamma$ increases, the intrinsic period will decrease. 
For large $\gamma$, this intrinsic length scales like $(1/\gamma)^{1/3}$ ({\it cf.} \cite{OK, C, Sp, ACO}) and hence as a very rough measure, we found it useful to set the initial domain size as 
\begin{equation}\label{domain-size}
L_\gamma \,= \, L_0\, \left( \frac{2}{\gamma}\right)^{1/3}.
\end{equation}
However, we want to ensure that we have reached 
a true global minimizer, and that boundary effects resulting from a ``bad" choice of domain size have not prevented us from reaching this goal. 
Thus at the end of each run,  we consider the field $u$ which has been numerically computed on $[0, L]^2$, rescale it to a unit domain $[0,1]^2$, and 
examine the rescaled form of the
energy per unit area:
\begin{align*}  \tilde{E} (L, \tilde{u},  \gamma, m)   \, & = \,  \frac{1}{L^2}\int_{[0,1]^2}  \frac{1}{\gamma^2} |\nabla \tilde{u} |^2  \, d {x} \,+\,
\int_{[0,1]^2}\frac{(1 - \tilde{u}^2)^2}{4} \, d {x}
 \, \\
& \quad + \, 
L^2 \int_{[0,1]^2}\int_{[0,1]^2} \, G({ x}, { y})  \left(\tilde{u}( x)  - {m}  \right)
\left(\tilde{u}( y)  - m \right) \,    d { x} \, d{y} \\
 & = \, \frac{1}{L^2\gamma^2} I_1(\tilde{u}) \, + \, I_2(\tilde{u}) \, + \, L^2 I_3(\tilde{u},m).  
\end{align*}
We then minimize $\tilde{E}$ with respect to $L$ for fixed $\tilde{u}$, $\gamma$ and $m$. 
With this choice of optimal $L$, we then integrate for a
small time further to ensure that the energy truly is lower.  This also gives us an additional
mechanism to find meta-stable states as they will have either a higher energy per unit volume or
drastically different optimal box size than their near neighbors in the $(m,\gamma)$ plane.

\subsection{Automatic pattern identification}\label{automatic}
To aid in the automatic classification of evolved states we examine the distribution of peaks 
with similar wave length to that of the dominant mode. In particular, we define and analyze 
the following function
\[g(\theta) = \int \exp{-\frac{|(k_x,k_y)- (\sin \theta, \cos \theta) k^*|^2}{\varepsilon}} \hat{v}(k)\, dk.\]
If $ g(\theta) $ (where $\theta = \tan^{-1} \frac{k_x}{k_y}$ is the polar angle) is found to have two dominant peaks for $0 \le \theta < 2 \pi$, the state is classified as a stripe.  If it has six then it is classified as a hexagonally packed circular state.   By symmetry, these are the $n=1$ (lamellae) and 3 (hex packed spots) cases from Section 2 respectively. 
Labyrinthine structures do not display either of these patterns.
This method is not perfect, working only in about $95\%$ of cases.

\subsection{Summary of algorithm}

In summary, we integrate (\ref{pde}) using the following algorithm:

\noindent
\begin{enumerate}
	\item Set $L = \left(\frac{2}{\gamma}\right)^{1/3}12\pi$.
	\item Choose random initial data in $[-1-m,1-m]$ with mean $m$.
\item Integrate (\ref{pde}) using ETDRK4 until $t = t_1$.
\item Integrate (\ref{pde}) using ETDRK4 with spectral weighting for $t_1 < t < t_2$.
\item Integrate (\ref{pde}) using ETDRK4 with added white noise for $t_2 < t < t_3$.
\item Integrate (\ref{pde}) using the nonlinear gradient stable scheme for $t_3 < t < t_4$.
\item Domain size variation/selection of Section \ref{domainsize}.
\item Integrate (\ref{pde}) using the nonlinear gradient stable scheme for $t_4 < t < t_5$.
\item Classify the final state.
\end{enumerate}

This procedure may appear overly complicated but recall that we are attempting to minimize a  non-convex, non-local energy 
in a manner which does not use an
anticipated solution structure (as this is not possible in 3D), and in a way which is efficient,
reliable and automatic.  To do this we have combined two different time-stepping strategies with
two different annealing mechanisms.

\section{Unusual Agreement with Asymptotic Analysis}\label{agree}

The calculations of the global stability curves of  Section \ref{asy} give surprisingly good
 agreement with the results of our spectrally weighted numerical algorithm for
a larger range of $\gamma$ than might be initially expected. This agreement extends far beyond the immediate neighborhood of the order-disorder curve. Naturally, this should raise a certain amount of skepticism and concern that the weighted numerical method biases the solutions towards the linear dynamics. 
Let us address these concerns.  First, we note that the spectral weighting is only employed for part of run. For those parts, it is certainly possible that a certain bias is introduced towards the linear dynamics.  
However, recall that our fundamental goal here is the phase diagram of (\ref{OK}): that is in parameter space $\gamma$ vs 
$m$, characterizing the geometric morphology of the global minimizer. 
Even without the spectral weighting,  all runs gave final states in which the basic pattern dichotomy (spots vs stripes) was recognizable. 
The  spectral weighting  allowed  for final states with  {\it lower energy} and which were identical to the {\it desired} morphology (eg. 
straight stripes vs corrugated).   
 At certain stages of the evolution, it does introduce a bias towards lowering the energy and  if this bias  favors the ``linear regime" of (\ref{pde}), so be it.  In fact, all of these conclusions lead to the following observation: 
For the purposes of describing the basic geometric morphology of the ground state of (\ref{OK}), the linear behavior of  (\ref{pde}) gives very accurate predictions in a (finite) neighborhood of the order-disorder transition. This is particularly true for initial data which is a small perturbation of the 
homogeneous state.  That is, the dynamics are initially driven by the linear terms and then there is a large
energy barrier to overcome to transition from one phase to the other.  In fact, starting with small initial
data gives convergence to states as predicted by the asymptotics for $\gamma$ as large as 25.  Starting with
initial conditions with {\it large random deviations} generates the diagram presented here as it is with these less biased initial conditions that we our able to access the lowest energy configurations. The domination of the linear terms
in constructing stationary solutions to higher-order PDE is well known in the folklore but little studied
beyond formal observations \cite{Pel46}

However, there is another more fundamental reason why this problem is unusual:  The energy has a term
which keeps solutions in the linear regime.  Solutions to the standard Cahn-Hilliard equation immediately tend to values of $\pm 1 + {\mathcal O}(\varepsilon)$ where $\varepsilon$ is the interfacial thickness. In (\ref{OK}),  the long-range interactions prefer to keep solutions close to oscillations about the mean.  
But this is precisely the asymptotic Ansatz,  and thus it is no surprise that it does well to  capture the
solutions when they are still in this regime.  Table \ref{Table:Fluc} presents the magnitude of
fluctuation about the mean as $\gamma$ varies and shows there is a large range of $\gamma$ where
the solution is in the asymptotic regime.

\begin{table}
	\begin{center}
		\begin{tabular}{c|c|c}
		$\gamma$ & $m^*$ & $\frac{\mbox{mean}_m(\max(u) - \min(u))}{2}$\\
		\hline
		2.001&  0.0050 & 0.013\\
		2.01& 0.0157 & 0.043\\		
		2.1& .0485 & 0.109\\
		2.25&  0.0741 & 0.163\\
		2.5& 0.2582 & 0.240\\
		3& 0.3333 & 0.317 \\
		3.5& 0.3780 & 0.413 \\
		5& 0.4472 & 0.626\\
		10& 0.5164 & 0.956\\
		20& 0.5477 & 1.021\\
		\end{tabular}
	\end{center}
	\medskip
	\caption{Variations about the mean as $m^*$ varies along the energy minimizing branches (see Figure 
	\ref{Fig:ContComp}).}
	\label{Table:Fluc}
\end{table}

There is also a region of the phase diagram where the agreement is quite poor: for large $\gamma$ and
large $m$ there is a region where there are very well separated hex packed spots and a very small 
region where there are Cartesian packed spots as well.  Examples of both of these with $\gamma = 20$ are
 presented in Figure \ref{fig:WeirdSpots}. Neither of these behaviors appears for small $m^*$ and 
the Cartesian packed spots are predicted to never be local minimizers. However, SCFT does predict a
region of ``close-packed" spheres which would be the 3D analogue of this behavior.

\begin{figure}[!t]
\centerline{{\includegraphics[width=0.5\textwidth]{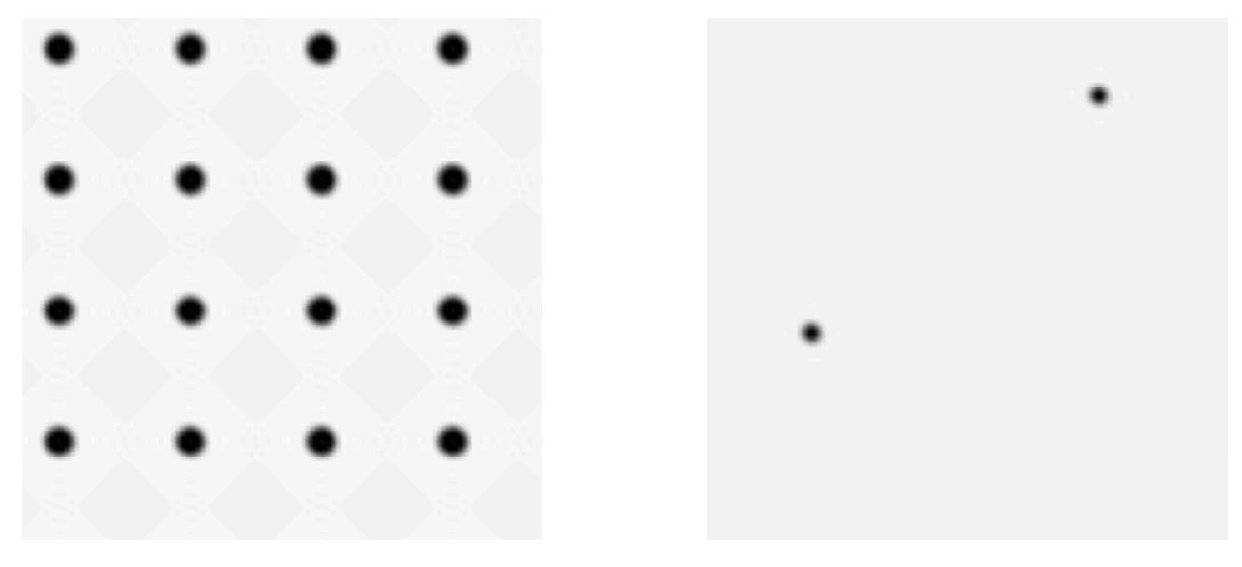}}}
\caption{Unusual solutions not predicted by linear analysis.  Both these cases have energy
lower than the mixed state but are not predicted by the linear analysis.}
\label{fig:WeirdSpots}
\end{figure}

\section{The Outlook for  3D}
The 3D problem is  considerably more complicated because of the multitude of both local and global minimizers. 
Preliminary experiments  suggest that our numerical method is reasonably successful throughout much of the phase plane.  
Phase boundaries on the other hand are much more difficult to capture and it is thus very encouraging that, at least in 2D, asymptotic analysis provides such a good indication of where these are close to the oder-disorder transition. Naturally, calculations for double gyroid and perforated lamellar symmetries are more complicated but are still tractable. 
As we saw in \cite{CPW} (though the simulations there were performed with a simpler, more ad hoc numerical method), metastability issues in 3D are even more complex and additional methods of 
{\it simulated annealing} may invariably be needed for a full phase diagram close to the order-disorder transition.  
Moreover our blind test for characterization of Section \ref{automatic} will need to be augmented with topological calculations on the phase boundary, e.g. the Euler characteristic, as was carried out in the recent work of Teramoto and Nishiura \cite{TN2}.

\medskip

{\bf Acknowledgments:}  We would like to thank one of the anonymous referees for their  questions and many suggestions which significantly improved the article. 
The authors were partially supported by the  
 NSERC (Canada) Discovery Grants program. 
A preliminary version of this work was contained in  the MSc thesis of Mirjana Maras at Simon Fraser University \cite{Ma}.

\end{document}